\documentclass{cocvnew}
\usepackage{eucal,enumerate,mathrsfs}
\usepackage{amsmath,amssymb,epsfig,bbm,amsthm}
\usepackage{ifthen}
\usepackage{pdfsync}
\newcommand{\HK}{\mathsf{H\kern-3pt K}}
\newcommand{\LET}{\mathsf{L\kern-2pt E\kern-2pt T}}
\newtheorem{ass}[thrm]{Assumption}
\numberwithin{equation}{section}
\usepackage{hyperref}
\begin{document}
%%-----------------------------
%%      the top matter
%%-----------------------------
\title{A Minimizing Movement approach to a class of scalar reaction-diffusion equations}%\thanks{...}\thanks{...}% At most 5 thanks
\runningtitle{A Minimizing Movement appraoch to reaction-diffusion equations}
\author{Florentine Catharina Flei$\ss$ner}\address{Technische Universit\"at M\"unchen email:
  \textsf{fleissne@ma.tum.de}.}
\begin{abstract}
The purpose of this paper is to introduce a Minimizing Movement approach to scalar reaction-diffusion equations of the form 
\begin{equation*}
\partial_t u \ = \ \Lambda\cdot \mathrm{div}[u(\nabla F'(u) + \nabla V)] \ - \ \Sigma\cdot (F'(u) + V) u, \quad \text{ in } (0, +\infty)\times\Omega,
\end{equation*}
with parameters $\Lambda, \Sigma > 0$ and no-flux boundary condition 
\begin{equation*}
u(\nabla F'(u) + \nabla V)\cdot {\sf n} \ = \ 0, \quad \text{ on } (0, +\infty)\times\partial\Omega,
\end{equation*}
which is built on their gradient-flow-like structure in the space $\mathcal{M}(\bar{\Omega})$ of finite nonnegative Radon measures on $\bar{\Omega}\subset\xR^d$, endowed with the recently introduced Hellinger-Kantorovich distance $\HK_{\Lambda, \Sigma}$. It is proved that, under natural general assumptions on $F: [0, +\infty)\to\xR$ and $V:\bar{\Omega}\to\xR$, the Minimizing Movement scheme 
\begin{equation*}
\mu_\tau^0:=u_0\mathscr{L}^d \in\mathcal{M}(\bar{\Omega}), \quad \mu_\tau^n \text{ is a minimizer for } \mathcal{E}(\cdot)+\frac{1}{2\tau}\HK_{\Lambda, \Sigma}(\cdot, \mu_\tau^{n-1})^2, \ n\in\xN, 
\end{equation*} 
%with discrete time step size $\tau > 0$ and 
for
\begin{equation*}
\mathcal{E}: \mathcal{M}(\bar{\Omega}) \to (-\infty, +\infty], \ \mathcal{E}(\mu):= \begin{cases} \int_\Omega{[F(u(x))+V(x)u(x)]\xdif x} &\text{ if } \mu=u\mathscr{L}^d, \\
+\infty &\text{ else},
\end{cases}
\end{equation*}
yields weak solutions to the above equation as the discrete time step size $\tau\downarrow 0$. Moreover, a superdifferentiability property of the Hellinger-Kantorovich distance $\HK_{\Lambda, \Sigma}$, which will play an important role in this context, is established in the general setting of a separable Hilbert space. 
\end{abstract}
%
%\begin{resume} ... \end{resume}
%
%\subjclass{0}
%
%\keywords{Gradient flows, Minimizing Movements, reaction-diffusion equations, Hellinger-Kantorovich distance, Optimal Transport}
%
\maketitle
%%-----------------------------
%%      your text
%%-----------------------------
\section*{Introduction}
At the beginning of the 90's, Ennio De Giorgi introduced the concept of \textit{Minimizing Movement} as ``natural meeting point" of many evolution problems from different reseach fields in mathematics \cite{DeGiorgi93}. He got his inspiration from the paper \cite{Almgren-Taylor-Wang93} by Almgren, Taylor and Wang. The concept involves the recursive minimization 
\begin{equation}\label{eq: MM scheme}
u_\tau^0:=u_0 \in \mathscr{S}, \quad u_\tau^n \text{ is a minimizer for } \Phi(\tau, u_\tau^{n-1}, \cdot), \ n\in\xN, 
\end{equation} 
of a given functional $\Phi: (0,1)\times \mathscr{S} \times \mathscr{S} \to [-\infty, +\infty]$ on a topological space $(\mathscr{S}, \sigma)$. The parameter $\tau > 0$ plays the role of discrete time step size. If a sequence $(u_\tau^n)_{n\in\xN}$ satisfies \eqref{eq: MM scheme}, we call the corresponding piecewise constant interpolation $u_\tau: [0, +\infty) \to \mathscr{S}, \ u_\tau(0)=u_0, \    
u_\tau(t) \equiv u_\tau^n \text{ for } t\in ((n-1)\tau, n\tau] \ (n\in\xN),$ a \textit{discrete solution}. The concept's purpose is to study the limit curves as $\tau\downarrow 0$.  

\begin{dfntn}[(Generalized) Minimizing Movement \cite{DeGiorgi93}]
A curve $u: [0, +\infty) \to \mathscr{S}$ is called Minimizing Movement for $\Phi$ with initial datum $u_0$ (short $u\in \mathrm{MM}(\Phi; u_0)$) if there exist discrete solutions $u_\tau$ to \eqref{eq: MM scheme} (for $\tau > 0$ in a right neighbourhood of $0$) such that $u_\tau^0 = u_0 = u(0)$ and $u_\tau(t)\stackrel{\sigma}{\rightharpoonup} u(t)$ for all $t > 0$ as $\tau\downarrow 0$. A curve $u: [0, +\infty) \to \mathscr{S}$ is called Generalized Minimizing Movement for $\Phi$ with initial datum $u_0$ (short $u\in \mathrm{GMM}(\Phi; u_0)$) if there exist a subsequence of time steps $(\tau_k)_{k\in\xN}, \ \tau_k\downarrow 0,$ and discrete solutions $u_{\tau_k}$ to \eqref{eq: MM scheme} such that $u_{\tau_k}^0 = u_0 = u(0)$ and $u_{\tau_k}(t)\stackrel{\sigma}{\rightharpoonup} u(t)$ for all $t > 0$ as $k\to\infty$.    
\end{dfntn}

\begin{xmpl}[Gradient flows in finite dimensional Euclidean space]\label{ex: GF in finite dimensions}
Let $\xQuaternion$ be a finite dimensional Euclidean space with norm $|\cdot|$ and $\mathcal{E}\in\xCone(\xQuaternion)$ satisfy the quadratic lower bound 
\begin{equation}
\exists A, B > 0: \quad \mathcal{E}(x) \geq -A-B|x|^2 \quad \text{ for all } x\in\xQuaternion.
\end{equation}
We apply the Minimizing Movement scheme \eqref{eq: MM scheme} to $\Phi(\tau, v, x):= \mathcal{E}(x) + \frac{1}{2\tau}|x-v|^2$. 
The necessary condition of first order leads to a discrete version of the classical gradient flow equation
\begin{equation}\label{eq: GF equation}
u'(t) = -\nabla\mathcal{E}(u(t)), \quad t\geq 0,
\end{equation}
and indeed, it is not difficult to see that every $u\in\mathrm{GMM}(\Phi; u_0)$ (which is a nonempty set) is a solution to \eqref{eq: GF equation} with initial datum $u_0\in\xQuaternion$. The statement also holds good if we replace $\mathcal{E}$ by $\mathcal{E}_\tau$ in $\Phi$, with $\mathcal{E}_\tau: \xQuaternion \to \xR$ converging to $\mathcal{E}$ in the Lipschitz semi-norm as $\tau\downarrow 0$.  Conversely, for every solution $u\in\xCone([0, +\infty);\xQuaternion)$ to \eqref{eq: GF equation} there exist functions $\mathcal{E}_\tau: \xQuaternion \to \xR$ $(\tau > 0)$ such that $\mathrm{Lip}[\mathcal{E}_\tau - \mathcal{E}] \to 0$ as $\tau\downarrow 0$ and $\mathrm{MM}(\Phi; u(0)) = \{u\} = \mathrm{GMM}(\Phi; u(0))$ for $\Phi(\tau, v, x):=\mathcal{E}_\tau(x) + \frac{1}{2\tau}|x-v|^2$, see \cite{sf2017}. This gives a full characterization of solutions to \eqref{eq: GF equation} as (Generalized) Minimizing Movements.    
\end{xmpl}
De Giorgi's concept of Minimizing Movements has a wide range of applications in analysis, geometry, physics and numerical analysis, and we refer to \cite{Almgren-Taylor-Wang93, DeGiorgi93, braides2012local, mielke2011differential, Mielke-Rindler09, mielke2016balanced} for more examples. In this paper the focus will be on the Minimizing Movement approach to gradient flows. \\

The Minimizing Movement scheme from Example \ref{ex: GF in finite dimensions} can be adapted for a general metric and non-smooth setting: 
Let a functional $\mathcal{E}: \mathscr{S}\to (-\infty, +\infty]$ on a complete metric space $(\mathscr{S}, d)$ be given and apply \eqref{eq: MM scheme} to 
\begin{equation}\label{eq: Phi}
\Phi(\tau, v, x):= \mathcal{E}(x) + \frac{1}{2\tau} d(x,v)^2. 
\end{equation}
It is proved in \cite{AmbrosioGigliSavare05} that, under natural coercivity assumptions, the set $\mathrm{GMM}(\Phi; u_0)$ is nonempty for every initial datum $u_0\in\{\mathcal{E} < +\infty\}$ and the Generalized Minimizing Movements are locally absolutely continuous curves satisfying the energy dissipation inequality 
\begin{equation}\label{eq: energy inequality}
\mathcal{E}(u(0)) - \mathcal{E}(u(t)) \ \geq \ \frac{1}{2} \int_0^t{|\partial^-\mathcal{E}|(u(r))^2 \mathrm d r} + \frac{1}{2}\int_0^t{|u'|(r)^2 \mathrm d r}
\end{equation}
for all $t > 0$, with $|\partial^-\mathcal{E}|$ denoting the relaxed slope of $\mathcal{E}$ (which can be viewed as a weak counterpart of the modulus of the gradient) and $|u'|$ the metric derivative of $u$ (see (\cite{AmbrosioGigliSavare05}, Chaps. 1 and 2) for the corresponding definitions or Sect. \ref{subsec: 3.3} in this paper for a brief overview). Under the additional assumption that the relaxed slope satisfies a kind of metric chain rule, equality can be proved in \eqref{eq: energy inequality}, see \cite{AmbrosioGigliSavare05}. The characterization of curves via such energy dissipation (in)equality corresponds with the notion of gradient flows in metric spaces which goes back to \cite{DeGiorgi-Marino-Tosques80, Degiovanni-Marino-Tosques85, Marino-Saccon-Tosques89}. It is equivalent to \eqref{eq: GF equation} if $\mathcal{E}\in\xCone(\xQuaternion)$ and $\xQuaternion$ is a finite dimensional Euclidean space. We refer to \cite{fleissner2016gamma} for further developments of the theory. 

Example \ref{ex: GF in finite dimensions} and the results in metric spaces motivate us to study a Minimizing Movement approach whenever gradient-flow-like structures are discovered and justify the interpretation of dynamics governed by an evolution equation as a gradient flow (for an energy functional on a metric space) if the corresponding (Generalized) Minimizing Movements are solutions to the evolution equation. 

The following Minimizing Movement approach to scalar diffusion equations of the form
\begin{equation}\label{eq: diffusion equation}
\partial_t u(t,x) \ = \ \mathrm{div}[u(t,x)(\nabla F'(u(t,x)) + \nabla V(x))], \quad t > 0, \ x\in\xR^d, 
\end{equation}  
has its origin in the papers \cite{JordanKinderlehrerOtto97, JordanKinderlehrerOtto98} by Jordan, Kinderlehrer and Otto, was examined by Ambrosio, Gigli and Savar\'e in \cite{AmbrosioGigliSavare05} and has been taken in many applications (see e.g. [\cite{AmbrosioGigliSavare05}, Chapt. 11] and the references therein): The space $\mathcal P_2(\xR^d)$ of Borel probability measures with finite second order moments (i.e. $\int_{\xR^d}{|x|^2\xdif \mu} < +\infty$) is endowed with the quadratic Wasserstein distance $\mathcal W_2$, 
\begin{equation}
\mathcal W_2(\mu_1, \mu_2)^2 := \min_{\gamma\in P(\mu_1, \mu_2)} \int_{\xR^d\times\xR^d}{|x-y|^2\xdif \gamma}, \quad \mu_i\in\mathcal P_2(\xR^d), 
\end{equation}
with $P(\mu_1, \mu_2)$ being the set of Borel probability measures on $\xR^d\times\xR^d$ whose first and second marginals coincide with $\mu_1$ and $\mu_2$ respectively (see e.g. \cite{Villani03, Villani09} for a detailed account of the theory of Optimal Transport and Wasserstein distances). The functional $\mathcal{E}: \mathcal P_2(\xR^d) \to (-\infty, +\infty]$,  
\begin{equation}
\mathcal{E}(\mu):=
\begin{cases}
\int_{\xR^d}^{}{[F(u(x)) + V(x)u(x)]\xdif x} &\text{ if } \mu = u\mathscr L^d \ (\mathscr L^d \text{ d-dimensional Lebesgue measure}), \\
+\infty &\text{ else}, 
\end{cases} 
\end{equation}
is defined on $(\mathcal P_2(\xR^d), \mathcal W_2)$ and $\Phi$ is defined according to \eqref{eq: Phi}. Under suitable assumptions on $F: [0, +\infty) \to \xR$ and $V: \xR^d\to\xR$, the corresponding Minimizing Movement scheme (often referred to as `JKO-scheme' in the literature) yields weak solutions to \eqref{eq: diffusion equation} (cf. the exemplary proof for $F(u)=u \log u$ and nonnegative $V\in\xCinfty(\xR^d)$ in \cite{JordanKinderlehrerOtto98} and Chaps. 10.1, 10.4 and 11.1.3 in \cite{AmbrosioGigliSavare05}): Such setting typically includes the assumptions that $F$ is convex, continuous with $F(0)=0$, differentiable in $(0,+\infty)$, has superlinear growth and is bounded from below by $s\mapsto - C s^\lambda$ for some $\lambda > \frac{d}{d+2}, \ C > 0,$ as well as locally Lipschitz continuity and nonnegativity of $V$.  

This paper concerns scalar reaction-diffusion equations of the form
\begin{equation}\label{eq: reaction-diffusion equation}
\partial_t u(t,x) \ = \ \Lambda\cdot \mathrm{div}[u(t,x)(\nabla F'(u(t,x)) + \nabla V(x))] \ - \ \Sigma\cdot (F'(u(t,x)) + V(x)) u(t,x), \quad t > 0, \ x\in\Omega\subset\xR^d,
\end{equation}
with fixed parameters $\Lambda > 0$ and $\Sigma > 0$ (with which the diffusion part and the reaction part respectively are weighted). Moreover, the reaction part is governed by the growth / shrinkage rate $\mathfrak{G}(x, u):= -(F'(u)+V(x))$ (with $F: [0, +\infty)\to\xR, \ V:\Omega\to\xR$) which also affects the diffusion part since according to \eqref{eq: reaction-diffusion equation}, diffusion occurs along the gradient $\nabla \mathfrak{G}(x, u(x))$ from regions of lower to higher growth rate / from regions of higher to lower shrinkage rate. Equation \eqref{eq: reaction-diffusion equation} seems a likely model for describing the evolution in time of the density $u$ of some biological, ecological, economic, ... quantity in various cases in which there is not only diffusion but also generation and annihilation of mass and in which the motion of the particles, members of the species, ... is influenced by their tendency to move towards regions with the most ``favourable'' conditions (see e.g. (\cite{kondratyev2016new}, Sect. 4), \cite{kondratyev2016fitness, kondratyev2020nonlinear} and the references therein for applications). Reaction-diffusion equations of the above form are closely related to a distance on the space of finite nonnegative Radon measures which has been recently introduced independently of each other by three different teams \cite{liero2016optimal, liero2018optimal, kondratyev2016new, chizat2018interpolating, chizat2018unbalanced}. We follow the presentation of the distance by Liero, Mielke and Savar\'e \cite{liero2016optimal, liero2018optimal} who named it \textit{Hellinger-Kantorovich distance}.    
\subsection{The Hellinger-Kantorovich distance}
Let $(X, {\sf d})$ be a Polish space (i.e. a complete separable metric space) and let $\mathcal{M}(X)$ be the space of finite nonnegative Radon measures on it. The class of Hellinger-Kantorovich distances $\HK_{\Lambda, \Sigma} \ (\Lambda, \Sigma > 0)$ can be characterized by the Logarithmic Entropy-Transport problems
\begin{equation}\label{eq: LET}
\LET_{\Lambda, \Sigma}(\mu_1, \mu_2) := \min \Big\{\sum_{i=1}^{2}{\frac{4}{\Sigma}\int_X{(\sigma_i\log \sigma_i - \sigma_i + 1)\xdif \mu_i}} + \int_{X\times X}{{\sf c}_{\Lambda, \Sigma}({\sf d}(x_1, x_2))\xdif \gamma}: \ \gamma\in\mathcal{M}(X\times X), \ \gamma_i\ll \mu_i\Big\},  
\end{equation} 
for $\mu_i\in\mathcal{M}(X)$, 
with 
\begin{equation}
\sigma_i := \frac{\xdif \gamma_i}{\xdif\mu_i} \ (\gamma_i \text{ i-th marginal of } \gamma), \quad {\sf c}_{\Lambda, \Sigma}({\sf d}):= 
\begin{cases}
-\frac{8}{\Sigma}\log(\cos(\sqrt{\Sigma/(4\Lambda)} {\sf d})) &\text{ if } {\sf d} < \pi\sqrt{\Lambda/\Sigma},\\
+\infty &\text{ if } {\sf d} \geq \pi\sqrt{\Lambda/\Sigma}.
\end{cases}
\end{equation}
An optimal plan $\gamma$ (which exists by (\cite{liero2018optimal}, Thm. 3.3)) describes an optimal way of converting $\mu_1$ into $\mu_2$ (possibly having different total mass) by means of transport and creation / annihilation of mass, in view of the transportation cost function ${\sf c}_{\Lambda, \Sigma}({\sf d})$ and the entropy cost functions $\frac{4}{\Sigma}(\sigma_i\log \sigma_i - \sigma_i + 1)$. The bigger the parameter $\Lambda > 0$ is (for the same $\Sigma>0$), the more the system favours transport. The bigger $\Sigma > 0$ is (for the same $\Lambda > 0$), the more the system favours creation and annihilation of mass. 

\begin{prpstn}[cf. Cor. 7.14, Thms. 7.15, 7.17, 7.20, Lem. 7.8 in \cite{liero2018optimal} and Sect. 3 in \cite{liero2016optimal}]\label{prop: basics on HK}
For all $\Lambda, \ \Sigma > 0$: 
\begin{enumerate}[\rm (1)] 
\item $\HK_{\Lambda, \Sigma}: \mathcal{M}(X)\times \mathcal{M}(X) \to [0, +\infty), \ \HK_{\Lambda, \Sigma}(\mu_1, \mu_2) := \sqrt{\LET_{\Lambda, \Sigma}(\mu_1, \mu_2)}$, is a distance on $\mathcal{M}(X)$.
\item $\HK_{\Lambda, \Sigma}$ metrizes the weak topology on $\mathcal{M}(X)$ in duality with continuous and bounded functions $\phi: X\to\xR$ (short $\phi\in\xCzero_b(X))$, i.e. 
\begin{equation}\label{eq: HK metrizes weak topology}
\lim_{n\to\infty}\HK_{\Lambda, \Sigma}(\mu_n, \mu) = 0 \quad\text{ if and only if }\quad \lim_{n\to\infty}\int_X{\phi\xdif \mu_n} = \int_X{\phi\xdif \mu} \text{ for all } \phi\in\xCzero_b(X). 
\end{equation} 
\item $(\mathcal{M}(X), \HK_{\Lambda, \Sigma})$ is a complete metric space. 
\item Let $\eta_0$ denote the null measure. For all $\mu\in\mathcal{M}(X)$: 
\begin{equation}\label{eq: null measure}
\HK_{\Lambda, \Sigma}(\mu, \eta_0)^2 = \frac{4}{\Sigma}\mu(X).
\end{equation}  
\end{enumerate}
\end{prpstn}
If $X = \xR^d$ or $X$ is a compact, convex subset of $\xR^d$ and ${\sf d}$ is induced by the usual norm, then a representation formula \`a la Benamou-Brenier can be proved for $\HK_{\Lambda, \Sigma}$ (see (\cite{liero2018optimal}, Thms. 8.18, 8.20; \cite{liero2016optimal}, Thm. 3.6(v))): 
\begin{equation}\label{eq: rep}
\HK_{\Lambda, \Sigma}(\mu_1, \mu_2)^2 = \inf\Big\{\int_0^1{\int_X{(\Lambda|\nabla\xi(t,x)|^2 + \Sigma|\xi(t, x)|^2)\xdif \mu_t(x)}\xdif t}: \ \mu_1 \stackrel{(\mu, \xi)}{\rightsquigarrow}\mu_2 \Big\}
\end{equation}    
where $\mu_1\stackrel{(\mu, \xi)}{\rightsquigarrow}\mu_2$ means that $\mu: [0, 1]\to\mathcal{M}(X)$ is a continuous curve connecting $\mu(0)=\mu_1$ and $\mu(1)=\mu_2$ and satisfying the continuity equation with reaction $\partial_t\mu_t = -\Lambda\mathrm{div}(\mu_t\nabla\xi_t) + \Sigma\mu_t\xi_t$, governed by $\xi: (0,1)\times X \to \xR$ with $\xi(t, \cdot)$ Lipschitz continuous and bounded for all $t\in(0,1)$, in duality with $\xCinfty$-functions with compact support in $(0.1)\times X$, i.e.
\begin{equation}\label{eq: continuity equation with reaction}
\int_0^1{\int_X{(\partial_t\psi(t,x) + \Lambda\nabla\psi(t,x)\cdot\nabla\xi(t,x) + \Sigma \psi(t,x)\xi(t,x))\xdif \mu_t(x)}\xdif t} \ = \ 0 \quad \text{ for all }\psi\in\xCinfty_c((0,1)\times X). 
\end{equation}
Hence, on the set $\{\mu\in\mathcal{M}(X): \ \mu = u\mathscr L^d\}$ of absolutely continuous Radon measures with respect to the Lebesgue measure, $\HK_{\Lambda, \Sigma}$ can be identified with the dissipation distance $\mathcal{D}_{\mathbb{K}_{\Lambda, \Sigma}}$,
\begin{equation}\label{eq: rep2}
\mathcal{D}_{\mathbb{K}_{\Lambda, \Sigma}}(u_1, u_2)^2 := \inf\Big\{\int_0^1{<\xi(t), \mathbb{K}_{\Lambda, \Sigma}(u_t)\xi(t)> \xdif t} \ : \ \partial_t u_t = \mathbb{K}_{\Lambda, \Sigma}(u_t)\xi(t), \ u_1\stackrel{(u, \xi)}{\rightsquigarrow}u_2 \Big\},
\end{equation}
generated by the Onsager operator $\mathbb{K}_{\Lambda, \Sigma}(u)\xi := -\Lambda\mathrm{div}(u\nabla\xi) + \Sigma u\xi$, which suggests a gradient-flow-like structure of \eqref{eq: reaction-diffusion equation}  associated with the energy functional $\mathcal{E}: \mathcal{M}(\Omega) \to (-\infty, +\infty]$, 
\begin{equation}\label{eq: mathcal E}
\mathcal{E}(\mu):=
\begin{cases} \int_{\Omega}^{}{[F(u(x)) + V(x)u(x)]\xdif x} &\text{ if } \mu=u\mathscr L^d,\\
+\infty &\text{ else, } 
\end{cases}
\end{equation}
on $(\mathcal{M}(\Omega), \HK_{\Lambda, \Sigma})$
(for details we refer to Sect. 2 in \cite{liero2016optimal}, \cite{mielke2011gradient}, Sect. 3.2 in \cite{kondratyev2016new}, and Otto's Riemannian formalism for $(\mathcal{P}_2(\xR^d), \mathcal W_2)$ in \cite{Otto01}).  

To handle such equation (in a weak form), Gallou\"et and Monsaigeon proposed a `JKO splitting scheme' in \cite{gallouet2017jko}, in which one step $\mu_\tau^n \curvearrowright \mu_\tau^{n+1}$ consists of two substeps
\begin{eqnarray*}
\mu_\tau^{n+1/2} &\text{ is a minimizer for } \mathcal{E}(\cdot) + \frac{1}{2\tau}\mathcal{W}_2(\cdot, \mu_\tau^n)^2, \\
\mu_\tau^{n+1} &\text{ is a minimizer for } \mathcal{E}(\cdot) + \frac{1}{2\tau}\mathsf{He}(\cdot, \mu_\tau^{n+1/2})^2 
\end{eqnarray*}  
(for $\Lambda = \Sigma = 1$), and which is justified by the interpretation of the Hellinger-Kantorovich distance as infimal convolution of the Kantorovich-Wasserstein distance $\mathcal{W}_2$ and the Hellinger-Kakutani/Fisher-Rao distance $\mathsf{He}$ (cf. \eqref{eq: rep} and \eqref{eq: rep2}). In this paper, we will work directly with the Hellinger-Kantorovich distance $\HK_{\Lambda, \Sigma}$ and take the `natural' Minimizing Movement approach to \eqref{eq: reaction-diffusion equation}, associated with 
\begin{equation}\label{eq: Phi our MM}
\Phi(\tau, \mu, \nu):= \mathcal{E}(\nu) + \frac{1}{2\tau}\HK_{\Lambda, \Sigma}(\nu, \mu)^2. 
\end{equation} 
Before presenting our results, we would like to mention \cite{chizat2017tumor} in which such approach has been taken for a particular equation of Hele-Shaw type, which serves as a model for tumour growth. The considerations therein are based on the special structure of the corresponding energy functional 
\begin{equation*}
\mathcal{E}(\mu) :=
\begin{cases}
-c\mu(\Omega) &\text{ if } \mu=u\mathscr{L}^d \text{ and } u\leq 1, \\
+\infty &\text{ else }
\end{cases}
\end{equation*}
(for $c>0$) and do not overlap with our analysis. 

\subsection{Our Minimizing Movement approach}\label{subsec: our MM approach}

Let $\Omega$ be an open, bounded, convex subset of $\xR^d$ with $\xCone$-boundary $\partial\Omega$, and for $\Lambda, \Sigma > 0$, let the space $\mathcal{M}(\bar{\Omega})$ of finite nonnegative Radon measures on its closure be endowed with the Hellinger-Kantorovich distance $\HK_{\Lambda, \Sigma}$, i.e. set $X:=\bar{\Omega}$ and ${\sf d}(x_1, x_2):=|x_1-x_2|$ (induced by the usual norm $|\cdot|$ on $\xR^d$) in \eqref{eq: LET} and Prop. \ref{prop: basics on HK}. We apply the Minimizing Movement scheme \eqref{eq: MM scheme} to 
\begin{equation}\label{eq: our MM}
\Phi(\tau, \mu, \nu):= \mathcal{E}(\nu) + \frac{1}{2\tau}\HK_{\Lambda, \Sigma}(\nu, \mu)^2,  
\end{equation}
\begin{equation}\label{eq: our E}
\mathcal{E}: \mathcal{M}(\bar{\Omega}) \to (-\infty, +\infty], \quad \mathcal{E}:=\mathcal{F} + \mathcal{V},  \quad 
\mathcal{F}(\mu):= 
\begin{cases}
\int_\Omega{F(u(x)) \xdif x} &\text{ if } \mu = u\mathscr{L}^d \\
+\infty &\text{ else } 
\end{cases}
, \quad \mathcal{V}(\mu) := \int_{\bar{\Omega}}{V(x)\xdif \mu},  
\end{equation}  
where $\mu=u\mathscr{L}^d$ means that $u: \Omega\to[0, +\infty)$ is Borel measurable and $\int_{\bar{\Omega}}{\phi(x)\xdif \mu} = \int_{\Omega}{\phi(x)u(x)\xdif x}$ for all $\phi\in\xCzero_b(\bar{\Omega})$. We prove that, under natural general assumptions on $F: [0, +\infty) \to \xR$ and $V: \bar{\Omega} \to \xR$ (comparable to the typical assumptions in the case of diffusion equations \eqref{eq: diffusion equation}, see above), the corresponding sets of Generalized Minimizing Movements $\mathrm{GMM}(\Phi; \mu_0)$ (for initial data $\mu_0\in\{\mathcal{E} < +\infty\}$) are nonempty and for every $\mu\in\mathrm{GMM}(\Phi; \mu_0)$ there is $u: [0, +\infty) \times \Omega \to [0, +\infty)$ such that $\mu(t) = u(t)\mathscr{L}^d$ for all $t\geq 0$ and $u$ solves the scalar reaction-diffusion equation
\begin{equation}\label{eq: our eq}
\partial_t u(t,x) \ = \ \Lambda\cdot \mathrm{div}[u(t,x)(\nabla F'(u(t,x)) + \nabla V(x))] \ - \ \Sigma\cdot (F'(u(t,x)) + V(x)) u(t,x), \quad t > 0, \ x\in\Omega,
\end{equation}
with no-flux boundary condition 
\begin{equation}\label{eq: our bc}
u(t,x)(\nabla F'(u(t,x)) + \nabla V(x))\cdot {\sf n}(x) \ = \ 0, \quad t>0, \ x\in\partial\Omega,
\end{equation}
in a weak form, see Thm. \ref{thm: main thm}. Here, ${\sf n}$ denotes the outward pointing unit normal vector field along $\partial\Omega$. We discuss our assumptions on $F$ and $V$ in Sect. \ref{subsec: 3.1}; they are satisfied for example if $V$ is Lipschitz continuous and $F(u) = c_1 u\log u$ ($c_1 > 0$) or $F(u) = -c_1u^q + c_2 u^p \ (c_1\geq 0, \ c_2 > 0, \ p>1, \ q\in(0,1))$, see Ex. \ref{ex: 3.8}. The key to proving our result is that we are able to establish a subdifferentiability property of the opposite Hellinger-Kantorovich distance $-\HK_{\Lambda, \Sigma}$ along certain directions. We can identify, for $\mu, \nu_0 \in \mathcal{M}(\bar{\Omega})$ and curves $h\mapsto\nu_h\in\mathcal{M}(\bar{\Omega})$ of the form 
\begin{equation}\label{eq: def of nuh}
\nu_h := (I+hv)_{\#}(1+hR)^2\nu_0 
\end{equation}
(where $v: \bar{\Omega}\to\xR^d, \ R:\bar{\Omega}\to\xR$ are bounded and the support of $v$ lies in $\Omega$),
elements of the Fr\'echet subdifferentials of the mappings
\begin{equation}\label{eq: 24}
h\mapsto -\frac{1}{2}\HK_{\Lambda, \Sigma}(\nu_h, \mu)^2
\end{equation}  
at $h=0$ and, setting $v:=\frac{4\Lambda}{\Sigma}\nabla\phi, \ R:=2\phi$ in \eqref{eq: def of nuh}, we can link them to the difference
\begin{equation}
\int_{\bar{\Omega}}{\phi\xdif \mu} - \int_{\bar{\Omega}}{\phi\xdif \nu_o},
\end{equation}
for any $\xCtwo$-function $\phi: \Omega\to\xR$ with compact support in $\Omega$, see Sect. \ref{subsec: directional derivative of HK}. Thereby, the possibility of establishing a discrete weak version of \eqref{eq: our eq} for discrete solutions to \eqref{eq: MM scheme} opens up. Further crucial points in our proof will be compactness issues, the passage to the limit $\tau\downarrow 0$ in the discrete weak version of \eqref{eq: our eq} and the Neumann boundary condition \eqref{eq: our bc}. 

The analysis of the Fr\'echet subdifferentials of the mappings \eqref{eq: 24} in Sect. \ref{subsec: directional derivative of HK} seems of independent interest and will be carried out for general separable Hilbert spaces. 

The plan for the paper is as follows. In Sect. \ref{sec: 1}, an equivalent characterization of the Hellinger-Kantorovich distance $\HK_{\Lambda, \Sigma}$ will be given, which will be useful for our study of subdifferentiability properties of $-\HK_{\Lambda, \Sigma}$ carried out in Sect. \ref{subsec: directional derivative of HK}. In Sect. \ref{sec: 3}, our Minimizing Movement approach to reaction-diffusion equations with no-flux boundary condition will be established. Our assumptions on $F$ and $V$ will be discussed in Sect. \ref{subsec: 3.1} and the proof of our main result will be given in Sect. \ref{subsec: 3.2}. In Sect. \ref{subsec: 3.3}, we will make some comments and go into future developments.    

\section{Hellinger-Kantorovich distance and Optimal Transportation on the cone}\label{sec: 1}
Let $(X, {\sf d})$ be a Polish space and define the geometric cone $\mathfrak{C}$ on $X$ as the quotient space 
\begin{equation}
\mathfrak{C}:= X\times [0, +\infty)/\sim
\end{equation}
with 
\begin{equation}
(x_1, r_1) \sim (x_2, r_2) \quad \Leftrightarrow \quad r_1 = r_2 = 0 \text{ or } r_1=r_2, \ x_1 = x_2
\end{equation}
for $x_i\in X, \ r_i\in[0, +\infty)$. The vertex $\mathfrak{o}$ (for $r=0$) and $[x,r]$ (for $x\in X$ and $r>0$) denote the corresponding equivalence classes, i.e. $\mathfrak{C} = \{[x,r] \ | \ x\in X, \ r>0\}\cup\{\mathfrak{o}\}$. 

In (\cite{liero2016optimal}, Sect. 3) and (\cite{liero2018optimal}, Sect. 7), the Logarithmic Entropy-Transport problem \eqref{eq: LET} for $\Lambda, \Sigma > 0$ is translated into a problem of optimal transportation on the cone governed by ${\sf d}_{\mathfrak{C}, \Lambda, \Sigma}: \mathfrak C \times \mathfrak C \to [0, +\infty),$
\begin{equation}
{\sf d}_{\mathfrak{C}, \Lambda, \Sigma}([x_1,r_1],[x_2, r_2])^2:= 
\frac{4}{\Sigma}\Big(r_1^2 + r_2^2 - 2 r_1 r_2 \cos\Big(\Big(\sqrt{\Sigma/4\Lambda}\ {\sf d}(x_1, x_2)\Big)\wedge\pi\Big)
\end{equation} 
(where $\mathfrak{o}$ is identified with $[\bar{x}, 0]$ for some $\bar{x}\in X$) which is a distance on $\mathfrak{C}$. 
The space $\mathcal{M}_2(\mathfrak{C})$ of finite nonnegative Radon measures on the cone with finite second order moments, i.e. $\int_{\mathfrak{C}}{{\sf d}_{\mathfrak{C}, \Lambda, \Sigma}([x,r],\mathfrak{o})^2\xdif \alpha([x,r])} < +\infty$, is endowed with an extended quadratic Kantorovich-Wasserstein distance $\mathcal{W}_{\mathfrak{C}, \Lambda, \Sigma}$, 
\begin{equation}\label{eq: Wasserstein on the cone}
\mathcal{W}_{\mathfrak{C}, \Lambda, \Sigma}(\alpha_1, \alpha_2)^2 := 
\begin{cases}
\min\Big\{\int_{\mathfrak{C}\times\mathfrak{C}}{{\sf d}_{\mathfrak{C}, \Lambda, \Sigma}([x_1,r_1],[x_2,r_2])^2\xdif \beta} \ | \ \beta\in M(\alpha_1, \alpha_2)\Big\} &\text{ if } \alpha_1(\mathfrak{C}) = \alpha_2(\mathfrak{C}), \\
+\infty &\text{ else, }
\end{cases}
\end{equation}
with $M(\alpha_1, \alpha_2)$ being the set of finite nonnegative Radon measures on $\mathfrak{C}\times \mathfrak{C}$ whose first and second marginals coincide with $\alpha_1$ and $\alpha_2$. Every measure $\alpha\in\mathcal{M}_2(\mathfrak{C})$ on the cone is assigned a measure $\mathfrak{h}\alpha\in\mathcal{M}(X)$ on $X$, 
\begin{equation}\label{eq: homogenous marginal}
\mathfrak{h}\alpha := {\sf x}_{\#}({\sf r}^2\alpha), \quad ({\sf x}, {\sf r}): \mathfrak{C}\to X\times [0, +\infty), \ ({\sf x}, {\sf r})([x, r]) := (x,r) \text{ for } [x,r]\in\mathfrak{C}, \ r > 0, \ ({\sf x}, {\sf r})(\mathfrak{o}):=(\bar{x},0),  
\end{equation}
i.e. $\int_X{\phi(x)\xdif (\mathfrak{h}\alpha)} = \int_{\mathfrak{C}}{{\sf r}^2\phi({\sf x})\xdif \alpha}$ for all $\phi\in\xCzero_b(X)$.
Note that the mapping $\mathfrak{h}: \mathcal{M}_2(\mathfrak{C}) \to \mathcal{M}(X)$ is not injective. It is proved in \cite{liero2018optimal} (see Probl. 7.4, Thm. 7.6, Lem. 7.9, Thm. 7.20 therein) that 
\begin{eqnarray}
\HK_{\Lambda, \Sigma}(\mu_1, \mu_2)^2 &=& \min\Big\{\mathcal W_{\mathfrak{C}, \Lambda, \Sigma}(\alpha_1, \alpha_2)^2 \ \Big| \ \alpha_i\in\mathcal{M}_2(\mathfrak{C}), \  \mathfrak{h}\alpha_i = \mu_i, \ i=1,2\Big\} \label{eq: cone 1} \\
&=& \min\Big\{\mathcal W_{\mathfrak{C}, \Lambda, \Sigma}(\alpha_1, \alpha_2)^2 + \frac{4}{\Sigma}\sum_{i=1}^{2}{(\mu_i - \mathfrak{h}\alpha_i)(X)} \ \Big| \ \alpha_i\in\mathcal{M}_2(\mathfrak{C}), \ \mathfrak{h}\alpha_i \leq \mu_i, \ i=1,2 \Big\} \label{eq: cone},
\end{eqnarray}
and that every optimal plan $\gamma\in\mathcal{M}(X\times X)$ for the Logarithmic Entropy-Transport problem \eqref{eq: LET} gives rise to a solution to \eqref{eq: cone} and vice versa. Moreover, if $\beta\in\mathcal{M}(\mathfrak{C}\times\mathfrak{C})$ is a solution to the transport problem (\eqref{eq: cone 1}, \eqref{eq: Wasserstein on the cone}) (which exists by (\cite{liero2018optimal}, Thm. 7.6)) or it is a solution to (\eqref{eq: cone}, \eqref{eq: Wasserstein on the cone}), then 
\begin{equation}\label{eq: transport}
\beta\Big(\Big\{([x_1, r_1], [x_2,r_2])\in\mathfrak{C}\times\mathfrak{C}: \ r_1, r_2 > 0, \ {\sf d}(x_1, x_2) > \pi\sqrt{\Lambda/\Sigma}\Big\}\Big)= 0,
\end{equation}  
(see (\cite{liero2018optimal}, Lem. 7.19)). 

This equivalent characterization of the Hellinger-Kantorovich distance $\HK_{\Lambda, \Sigma}$ has proved extremely useful in \cite{liero2018optimal} for the examination of structural properties; for example, the proofs therein of Prop. \ref{prop: basics on HK} (1), (2), (3) are based on it. The fact that all the information on transport of mass and creation / annihilation of mass according to \eqref{eq: LET} lies in a pure transportation problem is also a good starting point for our analysis of subdifferentiability properties of $-\HK_{\Lambda, \Sigma}$ in Sect. \ref{subsec: directional derivative of HK}. In this context, geodesics in $(\mathfrak{C}, {\sf d}_{\mathfrak{C}, \Lambda, \Sigma})$ will play a certain role, i.e. curves $\eta: [0, 1]\to \mathfrak{C}$ satisfying ${\sf d}_{\mathfrak{C}, \Lambda, \Sigma}(\eta(t), \eta(s)) = |t-s| {\sf d}_{\mathfrak{C}, \Lambda, \Sigma}(\eta(0), \eta(1))$ for all $s, t\in[0,1]$. 

We construct a geodesic which connects $[x_1, r_1]$ and $[x_2, r_2]$, supposing that ${\sf d}(x_1, x_2) \leq \pi\sqrt{\Lambda/\Sigma}, \ r_1, r_2 > 0,$ and that there exists a geodesic $\mathcal{X}$ in $(X, {\sf d})$ between $x_1$ and $x_2$, $x_1\neq x_2$ (cf. Sect. 8.1 in \cite{liero2018optimal}): Let us try to find functions $\mathcal{R}: [0, 1]\to [0, +\infty)$ and $\theta: [0, 1]\to [0, 1]$ such that $\eta: [0,1]\to\mathfrak{C}, \ \eta(t):=[\mathcal{X}(\theta(t)), \mathcal{R}(t)],$ is such geodesic. We note that 
\begin{equation}\label{eq: z1}
{\sf d}_{\mathfrak{C}, \Lambda, \Sigma}(\eta(t), \eta(s))^2 \ = \ \frac{4}{\Sigma}\Big(\mathcal{R}(t)^2 + \mathcal{R}(s)^2 - 2\mathcal{R}(s)\mathcal{R}(t)\cos\Big(\sqrt{\Sigma/4\Lambda}\ |\theta(t) - \theta(s)| {\sf d}(x_1, x_2)\Big)\Big) \ = \ |z(t) - z(s)|^2_{\xC},
\end{equation}
where $z: [0, 1]\to\xC$ is the curve in the complex plane $\xC$ defined as 
\begin{equation}\label{eq: z}
z(t):=\frac{2}{\sqrt{\Sigma}}\mathcal{R}(t)\exp\Big(i\theta(t)\sqrt{\Sigma/4\Lambda} \ {\sf d}(x_1, x_2)\Big),
\end{equation}
and $|\cdot|_\xC$ denotes the absolute value for complex numbers. If $z$ is a geodesic in the complex plane between $z_1:=\frac{2}{\sqrt{\Sigma}} r_1$ and $z_2:=\frac{2}{\Sigma}r_2\exp\Big(i\sqrt{\Sigma/4\Lambda} \ {\sf d}(x_1, x_2)\Big)$, then, according to \eqref{eq: z1}, the corresponding curve $\eta$ in the cone space, $\eta(t):=[\mathcal{X}(\theta(t)), \mathcal{R}(t)]$, is a geodesic between $[x_1, r_1]$ and $[x_2, r_2]$. Hence, the condition
\begin{equation}
z(t) = z_1 + t(z_2-z_1) \quad \text{ for all } t\in[0,1],
\end{equation}
which means
\begin{eqnarray}
\mathcal{R}(t)\cos\Big(\theta(t)\sqrt{\Sigma/4\Lambda} \ {\sf d}(x_1, x_2)\Big) &=& r_1 + t\Big(r_2\cos\Big(\sqrt{\Sigma/4\Lambda} \ {\sf d}(x_1, x_2)\Big) - r_1\Big) \quad \text{ for all } t\in[0,1], \label{eq: choice for R and theta 1} \\
\mathcal{R}(t)\sin\Big(\theta(t)\sqrt{\Sigma/4\Lambda} \ {\sf d}(x_1, x_2)\Big) &=& tr_2\sin\Big(\sqrt{\Sigma/4\Lambda} \ {\sf d}(x_1, x_2)\Big) \quad \text{ for all } t\in[0,1], \label{eq: choice for R and theta}
\end{eqnarray}  
yields an appropriate choice for $\mathcal{R}: [0, 1] \to [0, +\infty)$ and $\theta: [0,1]\to[0,1]$. It is not difficult to see that, by \eqref{eq: choice for R and theta 1} and \eqref{eq: choice for R and theta}, $\mathcal{R}$ and $\theta$ are smooth functions, their first derivatives satisfy
\begin{equation}\label{eq: first derivatives R and theta}
\frac{4}{\Sigma}(\mathcal{R}'(t))^2 + \frac{1}{\Lambda}\mathcal{R}(t)^2(\theta'(t))^2{\sf d}(x_1, x_2)^2 \ = \ {\sf d}_{\mathfrak{C}, \Lambda, \Sigma}([x_1, r_1], [x_2, r_2])^2 \quad \text{ for all } t\in(0,1),
\end{equation} 
and they are right differentiable at $t=0$ with right derivatives
\begin{equation}\label{eq: right derivative at 0}
\theta'_+(0) = \frac{r_2}{r_1}\frac{\sin(\sqrt{\Sigma/4\Lambda} \ {\sf d}(x_1, x_2))}{\sqrt{\Sigma/4\Lambda} \ {\sf d}(x_1, x_2)} \quad \text{ and } \quad \mathcal{R}'_+(0) = r_2 \cos\Big(\sqrt{\Sigma/4\Lambda} \ {\sf d}(x_1, x_2)\Big) - r_1. 
\end{equation}
Finally, we note that the curves $\eta: [0, 1]\to \mathfrak{C}$,
\begin{equation}
\eta(t):=
\begin{cases}
\mathfrak{o} &\text{ if } r_1 + t(r_2 - r_1) = 0, \\
 [x, r_1 + t(r_2 - r_1)] &\text{ else}, 
\end{cases}
\end{equation}
are geodesics in $(\mathfrak{C}, {\sf d}_{\mathfrak{C}, \Lambda, \Sigma})$, connecting $[x, r_1]$ (for $r_1>0$) or the vertex $\mathfrak{o}$ (for $r_1=0$)  with $[x, r_2]$ (for $r_2 > 0$) and with $\mathfrak{o}$ (for $r_2=0$). They take the above form $\eta(t)=[\mathcal{X}(\theta(t)), \mathcal{R}(t)]$ if we set $\mathcal{X} \equiv x$ (and identify the vertex with $[x, 0]$ if necessary),
\begin{equation}\label{eq: R and theta special cases}
\theta \equiv 0 \quad \text{ and } \quad \mathcal{R}(t):= r_1 + t(r_2 - r_1),
\end{equation}  also satisfying \eqref{eq: first derivatives R and theta} and the second part of \eqref{eq: right derivative at 0}.
 
\section{Superdifferentiability properties of the Hellinger-Kantorovich distance} \label{subsec: directional derivative of HK}  
Whenever a new distance is introduced, the question of differentiability properties arises. For the class of Hellinger-Kantorovich distances $\HK_{\Lambda, \Sigma}, \ \Lambda, \Sigma > 0,$ there has not been a corresponding analysis in the literature yet. In this section, we restrict ourselves to studying the superdifferentiability of $\HK_{\Lambda, \Sigma}$ (i.e. subdifferentiability of $-\HK_{\Lambda, \Sigma}$) along basic directions \eqref{eq: def of nuh} while we postpone studying the differentiability along general absolutely continuous curves to a subsequent paper. 

Let $\xQuaternion$ be a separable Hilbert space with scalar product $\langle\cdot, \cdot\rangle$ and norm $||\cdot||:=\sqrt{\langle\cdot, \cdot\rangle}$ and let $X\subset\xQuaternion$ be closed and convex. The couple $(X, {\sf d})$ with ${\sf d}(x_1, x_2) := ||x_1 - x_2||$ forms a Polish space. For $\Lambda, \Sigma > 0$, let the space $\mathcal{M}(X)$ of finite nonnegative Radon measures on $X$ be endowed with the distance $\HK_{\Lambda, \Sigma}$. 

We fix bounded Borel measurable functions $v: X\to \xQuaternion$ and $R: X\to\xR$, supposing that, for $h$ in a neighbourhood $\mathcal{N}$ around $0$, the function $I+hv: X\to\xQuaternion$ maps $X$ into $X$, where $I$ denotes the identity mapping $I: \xQuaternion \to \xQuaternion, \ I(x):=x$, and $1+hR(x) > 0$ for all $x\in X$ (which is satisfied whenever $|h|$ is small enough since $R$ is bounded). We define, for a given $\nu_0\in\mathcal{M}(X)$, the curve $\mathcal{N}\ni h\mapsto\nu_h\in\mathcal{M}(X)$ as  
\begin{equation}\label{eq: our direction}
\nu_h:=(I+hv)_{\#}(1+hR)^2\nu_0,
\end{equation}   
i.e. $\int_X{\phi(x)\xdif \nu_h} = \int_X{\phi(x+hv(x))(1+hR(x))^2\xdif \nu_0}$ for all $\phi\in\xCzero_b(X)$. Our goal is to identify elements of the Fr\'echet subdifferential of 
\begin{equation}
h\mapsto -\frac{1}{2}\HK_{\Lambda, \Sigma}(\nu_h, \mu)^2 \quad\quad (\mu\in\mathcal{M}(X))
\end{equation}
at $h=0$. A good strategy for this is to examine the subdifferentiability issue on the level of the optimal transportation problem on the associated cone $\mathfrak{C}$ first. We refer to Sect. \ref{sec: 1} for notation and details about the optimal transport problem on $\mathfrak{C}$. 

\begin{dfntn}[Fr\'echet subdifferential]
We say that $\varsigma\in\xR$ belongs to the Fr\'echet subdifferential of a mapping $\mathcal{N}\ni h\mapsto f(h)\in\xR$ at $h=0$ if and only if 
\begin{equation}
\liminf_{h\to 0}\frac{f(h)-f(0)-\varsigma h}{|h|} \ \geq \ 0.
\end{equation} 
\end{dfntn} 

\begin{lmm}\label{lmm: subdif of Wasserstein on the cone}
For a given $\alpha_0\in\mathcal{M}_2(\mathfrak{C})$, we define the curve $\mathcal{N}\ni h\mapsto\alpha_h\in\mathcal{M}_2(\mathfrak{C})$ as
\begin{equation}\label{eq: alphah}
\alpha_h:= ([x+hv(x), r(1+hR(x))])_{\#}\alpha_0,
\end{equation}
i.e. $\int_{\mathfrak{C}}{\varphi([x,r])\xdif\alpha_h} = \int_{\mathfrak{C}}{\varphi([x+hv(x), r(1+hR(x))])\xdif\alpha_0}$ for all $\varphi\in\xCzero_b(\mathfrak{C})$. Let $\alpha_\star\in\mathcal{M}_2(\mathfrak{C})$ be given with $\alpha_\star(\mathfrak{C}) = \alpha_0(\mathfrak{C})$ and let $\beta_{0,\star}\in M(\alpha_0, \alpha_\star)$ be optimal in the definition of $\mathcal{W}_{\mathfrak{C},\Lambda,\Sigma}(\alpha_0, \alpha_\star)^2$ according to \eqref{eq: Wasserstein on the cone}. We suppose that $\beta_{0, \star}$ satisfies \eqref{eq: transport}. Then the Fr\'echet subdifferential of the mapping
\begin{equation}
h\mapsto -\frac{1}{2}\mathcal{W}_{\mathfrak{C}, \Lambda, \Sigma}(\alpha_h, \alpha_\star)^2
\end{equation}
at $h=0$ is nonempty and 
\begin{equation}\label{eq: subdif of Wasserstein on the cone}
\frac{4}{\Sigma}\int_{\mathfrak{C}\times\mathfrak{C}}{\Big[-r_1^2R(x_1)+r_1r_2R(x_1)\cos(\sqrt{\Sigma/4\Lambda}||x_1-x_2||)+r_1r_2\sqrt{\Sigma/4\Lambda} \ \langle S_{\Lambda,\Sigma}(x_1,x_2), v(x_1)\rangle\Big]\xdif\beta_{0,\star}}
\end{equation}
belongs to it, where
\begin{equation}\label{eq: SLambdaSigma}
S_{\Lambda, \Sigma}(x_1, x_2):= \begin{cases}
\frac{\sin(\sqrt{\Sigma/4\Lambda}||x_1-x_2||)}{||x_1-x_2||}(x_2-x_1) &\text{ if } x_1\neq x_2, \\
0 &\text{ if } x_1=x_2.
\end{cases}  
\end{equation}  
\end{lmm}
\begin{proof}
First of all, we note that $\alpha_h(\mathfrak{C}) = \alpha_0(\mathfrak{C})$ and we define 
\begin{equation*}
\beta_{h, \star}:=([x_1+hv(x_1), r_1(1+hR(x_1))], [x_2, r_2])_{\#}\beta_{0, \star},
\end{equation*}
i.e. $\int_{\mathfrak{C}\times\mathfrak{C}}{\varphi([x_1, r_1],[x_2,r_2])\xdif\beta_{h,\star}} = \int_{\mathfrak{C}\times\mathfrak{C}}{\varphi([x_1+hv(x_1), r_1(1+hR(x_1))], [x_2, r_2])\xdif\beta_{0, \star}}$ for all $\varphi\in\xCzero_b(\mathfrak{C}\times\mathfrak{C})$. Since $\beta_{h,\star}\in M(\alpha_h, \alpha_\star)$, we have 
\begin{eqnarray*}
\mathcal{W}_{\mathfrak{C}, \Lambda, \Sigma}(\alpha_h, \alpha_\star)^2 &\leq& \int_{\mathfrak{C}\times\mathfrak{C}}{{\sf d}_{\mathfrak{C}, \Lambda, \Sigma}([x_1, r_1], [x_2, r_2])^2\xdif\beta_{h, \star}} \\ &=& \int_{\mathfrak{C}\times\mathfrak{C}}{{\sf d}_{\mathfrak{C}, \Lambda, \Sigma}([x_1+hv(x_1), r_1(1+hR(x_1))], [x_2, r_2])^2\xdif\beta_{0, \star}},
\end{eqnarray*}
leading to 
\begin{eqnarray*}
&& \frac{1}{2}\Big(\mathcal{W}_{\mathfrak{C}, \Lambda, \Sigma}(\alpha_0, \alpha_\star)^2 - \mathcal{W}_{\mathfrak{C}, \Lambda, \Sigma}(\alpha_h, \alpha_\star)^2\Big) \\ &\geq& 
\frac{1}{2}\int_{\mathfrak{C}\times\mathfrak{C}}{\Big[{\sf d}_{\mathfrak{C}, \Lambda, \Sigma}([x_1, r_1], [x_2,r_2])^2 - {\sf d}_{\mathfrak{C}, \Lambda, \Sigma}([x_1+hv(x_1), r_1(1+hR(x_1))], [x_2, r_2])^2\Big]\xdif \beta_{0, \star}} \\
&=& \frac{4}{\Sigma}\int_{\mathfrak{C}\times\mathfrak{C}}{\Big[-r_1^2hR(x_1) - \frac{1}{2}h^2R(x_1)^2r_1^2 + r_1r_2hR(x_1)\cos((\sqrt{\Sigma/4\Lambda}||x_1 + hv(x_1) - x_2||)\wedge\pi) \Big]\xdif\beta_{0, \star}} \\
&+& \frac{4}{\Sigma}\int_{\mathfrak{C}\times\mathfrak{C}}{\Big[r_1r_2\cos((\sqrt{\Sigma/4\Lambda}||x_1 + hv(x_1) - x_2||)\wedge\pi)-r_1r_2\cos(\sqrt{\Sigma/4\Lambda}||x_1-x_2||)\Big]\xdif\beta_{0, \star}}
\end{eqnarray*}
By dominated convergence theorem, it follows that
\begin{eqnarray*}
&&\liminf_{h\downarrow0}\frac{\frac{1}{2}\mathcal{W}_{\mathfrak{C}, \Lambda, \Sigma}(\alpha_0, \alpha_\star)^2 - \frac{1}{2}\mathcal{W}_{\mathfrak{C}, \Lambda, \Sigma}(\alpha_h, \alpha_\star)^2}{h} \ \geq \ \\  &&\frac{4}{\Sigma}\int_{\mathfrak{C}\times\mathfrak{C}}{\Big[-r_1^2R(x_1)+r_1r_2R(x_1)\cos(\sqrt{\Sigma/4\Lambda}||x_1-x_2||)+r_1r_2\sqrt{\Sigma/4\Lambda} \ \langle S_{\Lambda, \Sigma}(x_1, x_2), v(x_1)\rangle\Big]\xdif\beta_{0,\star}} \\
&&\geq \ \limsup_{h\uparrow0}\frac{\frac{1}{2}\mathcal{W}_{\mathfrak{C}, \Lambda, \Sigma}(\alpha_0, \alpha_\star)^2 - \frac{1}{2}\mathcal{W}_{\mathfrak{C}, \Lambda, \Sigma}(\alpha_h, \alpha_\star)^2}{h} \quad , 
\end{eqnarray*}
which completes the proof. 
\end{proof}
\begin{prpstn}\label{prop: subdif of HK}
For a given $\nu_0\in\mathcal{M}(X)$, we define the curve $\mathcal{N}\ni h\mapsto \nu_h\in\mathcal{M}(X)$ as in \eqref{eq: our direction}. Let $\mu\in\mathcal{M}(X)$ be given and let $\beta_{0, \star}\in\mathcal{M}(\mathfrak{C}\times\mathfrak{C})$ be optimal in the definition of $\HK_{\Lambda, \Sigma}(\nu_0, \mu)^2$ according to \eqref{eq: cone}, \eqref{eq: Wasserstein on the cone}, with first marginal $\alpha_0\in\mathcal{M}_2(\mathfrak{C}), \ \mathfrak{h}\alpha_0\leq \nu_0$, and second marginal $\alpha_\star\in\mathcal{M}_2(\mathfrak{C}), \ \mathfrak{h}\alpha_\star\leq \mu$. Then the Fr\'echet subdifferential of the mapping 
\begin{equation}\label{eq: 50}
h\mapsto - \frac{1}{2}\HK_{\Lambda, \Sigma}(\nu_h, \mu)^2
\end{equation}
at $h=0$ is nonempty and
\begin{equation}\label{eq: 51}
\mathfrak{F}_{0, \star, v, R} \ - \ \frac{4}{\Sigma}\int_{X}{R(x)\xdif(\nu_0-\mathfrak{h}\alpha_0)}
\end{equation}
belongs to it, where $\mathfrak{F}_{0, \star, v, R}$ is defined as
\begin{equation}\label{eq: 52} \frac{4}{\Sigma}\int_{\mathfrak{C}\times\mathfrak{C}}{\Big[-r_1^2R(x_1)+r_1r_2R(x_1)\cos(\sqrt{\Sigma/4\Lambda}||x_1-x_2||)+r_1r_2\sqrt{\Sigma/4\Lambda} \ \langle S_{\Lambda,\Sigma}(x_1,x_2), v(x_1)\rangle\Big]\xdif\beta_{0,\star}},  
\end{equation}
with $S_{\Lambda, \Sigma}$ as in \eqref{eq: SLambdaSigma}. 
\end{prpstn}
\begin{proof}
We define the curve $\mathcal{N}\ni h \mapsto \alpha_h\in \mathcal{M}_2(\mathfrak{C})$ as in \eqref{eq: alphah}. It holds that
\begin{eqnarray*}
\int_{X}{\phi(x)\xdif(\mathfrak{h}\alpha_h)} &=& \int_{\mathfrak{C}}{{\sf r}^2\phi({\sf x})\xdif\alpha_h}  = \int_{\mathfrak{C}}{{\sf r}^2(1+hR({\sf x}))^2\phi({\sf x}+hv({\sf x}))\xdif\alpha_0}  =  \int_{X}{(1+hR(x))^2\phi(x+hv(x))\xdif(\mathfrak{h}\alpha_0)} \\
&\leq& \int_{X}{(1+hR(x))^2\phi(x+hv(x))\xdif\nu_0} \ = \ \int_{X}{\phi(x)\xdif\nu_h}
\end{eqnarray*}
for all nonnegative bounded Borel functions $\phi: X\to \xR$, see \eqref{eq: homogenous marginal}, from which we infer that
\begin{equation*}
\mathfrak{h}\alpha_h \ \leq \ \nu_h \quad\text{ and }\quad (\nu_h - \mathfrak{h}\alpha_h)(X) \ = \ \int_X{(1+hR(x))^2\xdif(\nu_0-\mathfrak{h}\alpha_0)}. 
\end{equation*} 
Hence, we have
\begin{eqnarray*}
\frac{1}{2}\Big(\HK_{\Lambda, \Sigma}(\nu_0, \mu)^2-\HK_{\Lambda, \Sigma}(\nu_h, \mu)^2\Big)  &\geq&  \frac{1}{2}\Big(\mathcal{W}_{\mathfrak{C}, \Lambda, \Sigma}(\alpha_0, \alpha_\star)^2 - \mathcal{W}_{\mathfrak{C}, \Lambda, \Sigma}(\alpha_h, \alpha_\star)^2\Big) \\ &+& \frac{2}{\Sigma}\int_{X}{(-2hR(x)-h^2R(x)^2)\xdif(\nu_0-\mathfrak{h}\alpha_0)},  
\end{eqnarray*}
and we conclude by applying Lem. \ref{lmm: subdif of Wasserstein on the cone}. 
\end{proof}
This result can also be expressed in terms of the Logarithmic Entropy-Transport characterization \eqref{eq: LET} of the Hellinger-Kantorovich distance $\HK_{\Lambda, \Sigma}$. By Thm. 7.20 in \cite{liero2018optimal}, every optimal plan $\gamma\in\mathcal{M}(X\times X)$ for the Logarithmic Entropy-Transport problem \eqref{eq: LET} gives rise to a solution $\beta\in\mathcal{M}(\mathfrak{C}\times \mathfrak{C})$ to the optimal transportation problem (\eqref{eq: cone}, \eqref{eq: Wasserstein on the cone}) on the cone. Therefore, we obtain
\begin{crllr}\label{cor: subdif of HK}
Let $\mu, \nu_0, \nu_h\in\mathcal{M}(X)$ be given as in Prop. \ref{prop: subdif of HK} and let $\gamma\in\mathcal{M}(X\times X)$ be optimal in the definition of $\HK_{\Lambda, \Sigma}(\nu_0, \mu)^2$ according to \eqref{eq: LET}, with first marginal $\gamma_0\ll \nu_0$ and second marginal $\gamma_\star\ll \mu$. We suppose that 
\begin{equation}\label{eq: Lebesgue decomposition} 
\nu_0 = \rho_0\gamma_0 + \nu_0^\bot \quad\text{ and } \quad \mu=\rho_\star\gamma_\star+\mu^\bot
\end{equation}
for Borel functions $\rho_0, \rho_\star : X\to [0, +\infty)$ and nonnegative finite Radon measures $\nu_0^\bot, \mu^\bot \in\mathcal{M}(X), \ \nu_0^\bot \bot \gamma_0, \ \mu^\bot \bot \gamma_\star$, i.e. $\int_X{\phi(x)\xdif \nu_0} = \int_X{\rho_0(x)\phi(x)\xdif \gamma_0} + \int_X{\phi(x)\xdif \nu_0^\bot}$ for all $\phi\in\xCzero_b(X)$ and there exists a Borel set $B_0\subset X$ such that $\nu_0^\bot(B_0) = 0 = \gamma_0(X\setminus B_0)$; similarly for $(\mu, \rho_\star,\gamma_\star, \mu^\bot)$. Then
\begin{equation*}
\mathfrak{F}_{0, \star, v, R} \ - \ \frac{4}{\Sigma}\int_{X}{R(x)\xdif\nu_0^\bot}
\end{equation*}
belongs to the Fr\'echet subdifferential of \eqref{eq: 50} at $h=0$, with $\mathfrak{F}_{0, \star, v, R}$ defined as
\begin{eqnarray*} && \frac{4}{\Sigma}\int_{X\times X}{\Big[-\rho_0(x_1)R(x_1)+\sqrt{\rho_0(x_1)\rho_\star(x_2)}R(x_1)\cos(\sqrt{\Sigma/4\Lambda}||x_1-x_2||)\Big]\xdif \gamma} \\  &+& \frac{4}{\Sigma}\int_{X\times X}{\sqrt{\Sigma/4\Lambda \ \rho_0(x_1)\rho_\star(x_2)}\ \langle S_{\Lambda,\Sigma}(x_1,x_2), v(x_1)\rangle\xdif\gamma}.   
\end{eqnarray*} 
\end{crllr}
\begin{proof}
We define
\begin{equation*}
\beta_{0, \star} := ([x_1, \sqrt{\rho_0(x_1)}], [x_2, \sqrt{\rho_\star(x_2)}])_{\#}\gamma \ \in \mathcal{M}(\mathfrak{C}\times\mathfrak{C}),
\end{equation*}
i.e. $\int_{\mathfrak{C}\times\mathfrak{C}}{\varphi([x_1, r_1], [x_2, r_2])\xdif\beta_{0, \star}} = \int_{X\times X}{\varphi([x_1, \sqrt{\rho_0(x_1)}], [x_2, \sqrt{\rho_\star(x_2)}])\xdif\gamma}$ for all $\varphi\in\xCzero_b(\mathfrak{C}\times\mathfrak{C})$. 
According to Thm. 7.20 (iii) in \cite{liero2018optimal}, $\beta_{0, \star}$ is a solution to (\eqref{eq: cone}, \eqref{eq: Wasserstein on the cone}). By Prop. \ref{prop: subdif of HK}, the claim is proved.  
\end{proof}
We remark that, according to Thm. 1.115 in \cite{fonseca2007modern} or Lem. 2.3 in \cite{liero2018optimal}, such `Lebesgue decomposition' \eqref{eq: Lebesgue decomposition} always exists. 

Our analysis of the Fr\'echet subdifferentials of the mappings \eqref{eq: 50} at $h=0$ will form the basis for the general study of differentiability properties of the Hellinger-Kantorovich distance $\HK_{\Lambda, \Sigma}$ in a subsequent paper. For the purposes of our Minimizing Movement approach to \eqref{eq: our eq}, \eqref{eq: our bc}, the results from Prop. \ref{prop: subdif of HK} and Cor. \ref{cor: subdif of HK} will be sufficient. We conclude this section with a link between the Fr\'echet subdifferential of \eqref{eq: 50} at $h=0$ for $v:=\frac{4\Lambda}{\Sigma}\nabla\phi, \ R:=2\phi$, and the difference
\begin{equation*}
\int_{X}{\phi(x)\xdif\mu} - \int_X{\phi(x)\xdif\nu_0}. 
\end{equation*}
\begin{prpstn}\label{prop: connection of subdif with MM}
Let $\phi: \xQuaternion\to\xR$ be a twice continuously differentiable function whose differentials of first and second order at $x\in\xQuaternion$ are represented by the gradient $\nabla\phi(x)\in\xQuaternion$ and the Hessian $\nabla^2\phi(x): \xQuaternion\to\xQuaternion$ respectively. We suppose that
\begin{equation}\label{eq: 54}
C_\phi := \sup_{x\in X}(|\phi(x)|+||\nabla\phi(x)||+|||\nabla^2\phi(x)|||) < +\infty,
\end{equation}
with 
\begin{equation*}
|||\nabla^2\phi(x)|||:=\sup\{||\nabla^2\phi(x)(v)||: \ v\in\xQuaternion, ||v||\leq 1\}. 
\end{equation*}
For $\nu_0, \mu\in\mathcal{M}(X)$, let $\beta_{0, \star}\in\mathcal{M}(\mathfrak{C}\times\mathfrak{C})$ be optimal in the definition of $\HK_{\Lambda, \Sigma}(\nu_0, \mu)^2$ according to \eqref{eq: cone}, \eqref{eq: Wasserstein on the cone}, with first marginal $\alpha_0\in\mathcal{M}_2(\mathfrak{C}), \ \mathfrak{h}\alpha_0\leq \nu_0$, and second marginal $\alpha_\star\in\mathcal{M}_2(\mathfrak{C}), \ \mathfrak{h}\alpha_\star\leq \mu$. Then the following holds good:
\begin{equation}\label{eq: 55}
\Big|\frac{4}{\Sigma}\Big(\int_X{\phi(x)\xdif\mu} - \int_X{\phi(x)\xdif\nu_0}\Big) - \Big(\mathfrak{F}_{0, \star, \phi} - \frac{8}{\Sigma}\int_X{\phi(x)\xdif(\nu_0-\mathfrak{h}\alpha_0)}\Big)\Big| \ \leq \ C_\phi(6+16\Lambda/\Sigma)\HK_{\Lambda, \Sigma}(\nu_0, \mu)^2,
\end{equation}
where $\mathfrak{F}_{0, \star, \phi}$ is defined as
\begin{equation}
\frac{4}{\Sigma}\int_{\mathfrak{C}\times\mathfrak{C}}{\Big[-2r_1^2\phi(x_1)+2r_1r_2\phi(x_1)\cos(\sqrt{\Sigma/4\Lambda}||x_1-x_2||)+r_1r_2\sqrt{4\Lambda/\Sigma}\ \langle S_{\Lambda,\Sigma}(x_1,x_2), \nabla\phi(x_1)\rangle\Big]\xdif\beta_{0,\star}} 
\end{equation}
(with $S_{\Lambda, \Sigma}$ as in \eqref{eq: SLambdaSigma}). 
\end{prpstn}
\begin{proof}
First of all, we note that
\begin{eqnarray*}
\int_X{\phi(x)\xdif\mu} - \int_X{\phi(x)\xdif\nu_0} \ = \  \int_{\mathfrak{C}\times\mathfrak{C}}{\Big[\phi(x_2)r_2^2 - \phi(x_1)r_1^2\Big]\xdif\beta_{0, \star}} \ + \int_X{\phi(x)\xdif(\mu-\mathfrak{h}\alpha_\star)} - \int_X{\phi(x)\xdif(\nu_0-\mathfrak{h}\alpha_0)}
\end{eqnarray*}
and 
\begin{equation*}
\Big|\frac{4}{\Sigma}\int_X{\phi(x)\xdif(\mu-\mathfrak{h}\alpha_\star)} - \frac{4}{\Sigma}\int_X{\phi(x)\xdif(\nu_0-\mathfrak{h}\alpha_0)} + \frac{8}{\Sigma}\int_X{\phi(x)\xdif(\nu_0-\mathfrak{h}\alpha_0)}\Big| \leq  C_\phi \HK_{\Lambda, \Sigma}(\nu_0, \mu)^2. 
\end{equation*}
Hence, all that remains is to find a suitable estimate of
$\Big|\frac{4}{\Sigma}\int_{\mathfrak{C}\times\mathfrak{C}}{\Big[\phi(x_2)r_2^2 - \phi(x_1)r_1^2\Big]\xdif\beta_{0, \star}} - \mathfrak{F}_{0, \star, \phi}\Big|$. We fix $([x_1, r_1], [x_2, r_2]) \in \mathfrak{C}\times \mathfrak{C}\setminus\{(\mathfrak{o}, \mathfrak{o})\}$ with $||x_1-x_2||\leq \pi \sqrt{\Lambda/\Sigma}$. Let $\eta: [0, 1]\to\mathfrak{C}, \ \eta(t):=[x_1+\theta(t)(x_2-x_1), \ \mathcal{R}(t)],$ be the geodesic between $[x_1, r_1]$ and $[x_2, r_2]$ in $(\mathfrak{C}, {\sf d}_{\mathfrak{C}, \Lambda, \Sigma})$, defined according to \eqref{eq: z1}-\eqref{eq: R and theta special cases}. Then the mapping $t\mapsto\chi(t):=\phi(x_1+\theta(t)(x_2-x_1))\mathcal{R}(t)^2$ is twice continuously differentiable with 
\begin{eqnarray*}
\chi'(t) = 2\mathcal{R}(t)\mathcal{R}'(t)\phi(x_1+\theta(t)(x_2-x_1))+\mathcal{R}(t)^2\theta'(t)\langle\nabla\phi(x_1+\theta(t)(x_2-x_1)),x_2-x_1\rangle \\
\chi''(t) = \Big(\frac{\mathrm{d}}{\mathrm{d}^2 t}\mathcal{R}(t)^2\Big)\phi(x_1+\theta(t)(x_2-x_1)) + 2\Big(\frac{\mathrm{d}}{\mathrm{d} t}\mathcal{R}(t)^2\Big)\theta'(t)\langle\nabla\phi(x_1+\theta(t)(x_2-x_1)), x_2-x_1\rangle \\
+ \mathcal{R}(t)^2\theta''(t)\langle\nabla\phi(x_1+\theta(t)(x_2-x_1)), x_2-x_1\rangle \\ + \mathcal{R}(t)^2(\theta'(t))^2\langle x_2-x_1, \nabla^2\phi(x_1+\theta(t)(x_2-x_1))(x_2-x_1)\rangle 
\end{eqnarray*}  
for $t\in (0,1)$ and it is right differentiable at $t=0$ with right derivative
\begin{eqnarray*}
\chi'_+(0) \ = \ \lim_{t\downarrow 0}\chi'(t) &=&2r_1\mathcal{R}'_+(0)\phi(x_1) + r_1^2\theta'_+(0)\langle\nabla\phi(x_1), x_2-x_1\rangle \\
&=& -2r_1^2\phi(x_1)+2r_1r_2\phi(x_1)\cos(\sqrt{\Sigma/4\Lambda} \ ||x_1-x_2||) + r_1r_2\sqrt{4\Lambda/\Sigma} \ \langle S_{\Lambda, \Sigma}(x_1, x_2), \nabla\phi(x_1)\rangle.
\end{eqnarray*}
Approximating $\chi$ by the first Taylor polynomial at $t=0$ yields 
\begin{equation*}
|\phi(x_2)r_2^2-\phi(x_1)r_1^2-\chi'_+(0)| \ = \ |\chi(1) - \chi(0) - \chi'_+(0)| \ \leq \ \sup_{t\in (0,1)}|\chi''(t)|. 
\end{equation*}
So let us fix $t\in (0,1)$ and estimate $|\chi''(t)|$. For this, we need to play with the first and second derivatives of $\mathcal{R}$ and $\theta$. We recall \eqref{eq: first derivatives R and theta}, which says
\begin{equation*}
\frac{4}{\Sigma}(\mathcal{R}'(t))^2 + \frac{1}{\Lambda}\mathcal{R}(t)^2(\theta'(t))^2||x_1-x_2||^2 \ = \ {\sf d}_{\mathfrak{C}, \Lambda, \Sigma}([x_1, r_1], [x_2, r_2])^2. 
\end{equation*}
It is not difficult to see that \eqref{eq: choice for R and theta 1} and \eqref{eq: choice for R and theta} imply
\begin{equation}\label{eq: 57}
\frac{4}{\Sigma} (\mathcal{R}'(t))^2 + \frac{4}{\Sigma} \mathcal{R}(t)\mathcal{R}''(t) \ = \ \frac{2}{\Sigma} \ \frac{\mathrm{d}}{\mathrm{d}^2 t}\mathcal{R}(t)^2 \ = \ {\sf d}_{\mathfrak{C}, \Lambda, \Sigma}([x_1, r_1], [x_2, r_2])^2. 
\end{equation}
We infer from \eqref{eq: first derivatives R and theta} and \eqref{eq: 57} that
\begin{equation*}
\frac{1}{\Lambda}\mathcal{R}(t)(\theta'(t))^2||x_1-x_2||^2 \ = \ \frac{4}{\Sigma}\mathcal{R}''(t)
\end{equation*}
(since $\mathcal{R}(t)\neq 0$) and, by taking the derivative in \eqref{eq: 57}, that
\begin{equation*}
3\mathcal{R}'(t)\mathcal{R}''(t) + \mathcal{R}(t)\mathcal{R}'''(t) \ = \ 0. 
\end{equation*}
It follows that
\begin{eqnarray*}
\frac{1}{\Lambda}\mathcal{R}(t)\mathcal{R}'(t)(\theta'(t))^2||x_1-x_2||^2 & = & \frac{4}{\Sigma}\mathcal{R}'(t)\mathcal{R}''(t) \ = \  -\frac{4}{3\Sigma}\mathcal{R}(t)\mathcal{R}'''(t)\\ &=& -\frac{1}{3\Lambda}\mathcal{R}(t)\mathcal{R}'(t)(\theta'(t))^2||x_1-x_2||^2 - \frac{2}{3\Lambda}\mathcal{R}(t)^2\theta'(t)\theta''(t)||x_1-x_2||^2. 
\end{eqnarray*}
Supposing that $x_1\neq x_2$ and that $\theta'(t)\neq 0$, we obtain
\begin{equation}\label{eq: 58}
\Big|\theta''(t)\mathcal{R}(t)^2||x_1-x_2||\Big| \ \leq \ 2\Big|\mathcal{R}'(t)\mathcal{R}(t)\theta'(t)||x_1-x_2||\Big| \ \leq \ (\mathcal{R}'(t))^2 + \mathcal{R}(t)^2(\theta'(t))^2||x_1-x_2||^2.  
\end{equation}
We note that \eqref{eq: 58} also holds good if $x_1=x_2$ or if $\theta'(t)=0$, since, taking the first and second derivative in \eqref{eq: first derivatives R and theta}, we see that $\theta'(t)=0$ implies $\theta''(t)=0$ or $x_1=x_2$. Indeed, if $\theta'(t)=0$, we have
\begin{eqnarray*}
\mathcal{R}''(t) \ = \ 0 \ = \ \mathcal{R}'''(t) 
\end{eqnarray*} 
by the above considerations, and thus, 
\begin{equation*}
0 \ = \ \frac{\mathrm{d}}{\mathrm{d}^2 t}{\sf d}_{\mathfrak{C}, \Lambda, \Sigma}([x_1, r_1], [x_2, r_2])^2 \ = \ \frac{2}{\Lambda}\mathcal{R}(t)^2(\theta''(t))^2||x_1-x_2||^2. 
\end{equation*}
Finally, we obtain
\begin{eqnarray*}
|\chi''(t)| &\leq& C_\phi \Big(\frac{\mathrm{d}}{\mathrm{d}^2 t}\mathcal{R}(t)^2 + \underbrace{2 \Big|\Big(\frac{\mathrm{d}}{\mathrm{d} t}\mathcal{R}(t)^2\Big)\theta'(t)||x_1-x_2||\Big|}_{\leq 2(\mathcal{R}'(t))^2+2\mathcal{R}(t)^2(\theta'(t))^2||x_1-x_2||^2} + \Big|\mathcal{R}(t)^2\theta''(t)||x_1-x_2||\Big|+\mathcal{R}(t)^2(\theta'(t))^2||x_1-x_2||^2\Big) \\
&\leq& C_\phi(5/4 \ \Sigma + 4\Lambda){\sf d}_{\mathfrak{C}, \Lambda, \Sigma}([x_1, r_1], [x_2, r_2])^2, 
\end{eqnarray*}
by applying \eqref{eq: first derivatives R and theta}, \eqref{eq: 57} and \eqref{eq: 58}, and the fact that $\frac{\mathrm{d}}{\mathrm{d} t}\mathcal{R}(t)^2 = 2\mathcal{R}'(t)\mathcal{R}(t)$. 

Note that $\beta_{0, \star}$ satisfies \eqref{eq: transport}. All in all, it follows that
\begin{eqnarray*}
\Big|\frac{4}{\Sigma}\int_{\mathfrak{C}\times\mathfrak{C}}{\Big[\phi(x_2)r_2^2 - \phi(x_1)r_1^2\Big]\xdif\beta_{0, \star}} - \mathfrak{F}_{0, \star, \phi}\Big| &\leq& C_\phi (5+16\Lambda/\Sigma)\int_{\mathfrak{C}\times\mathfrak{C}}{{\sf d}_{\mathfrak{C}, \Lambda, \Sigma}([x_1, r_1], [x_2, r_2])^2\xdif \beta_{0, \star}} \\
&\leq& C_\phi (5+16\Lambda/\Sigma) \HK_{\Lambda, \Sigma}(\nu_0, \mu)^2.
\end{eqnarray*}
The proof of Prop. \ref{prop: connection of subdif with MM} is complete. 
\end{proof}
\begin{rmrk}\label{rem: connection of subdif with MM}
Let $\phi: \xQuaternion\to \xR$ satisfy the assumptions of Prop. \ref{prop: connection of subdif with MM}. In addition, we suppose that $\phi$ has compact support within the interior of $X$ so that, setting $v:=\frac{4\Lambda}{\Sigma}\nabla\phi$ and $R:=2\phi$, we can define the curve $\mathcal{N}\ni h\mapsto \nu_h\in\mathcal{M}(X)$ according to \eqref{eq: our direction}. Then
\begin{equation}
\mathfrak{F}_{0, \star, \phi} - \frac{8}{\Sigma}\int_X{\phi(x)\xdif(\nu_0 - \mathfrak{h}\alpha_0)},
\end{equation}
defined as in Prop. \ref{prop: connection of subdif with MM}, belongs to the Fr\'echet subdifferential of 
\begin{equation}
h\mapsto -\frac{1}{2}\HK_{\Lambda, \Sigma}(\nu_h, \mu)^2
\end{equation}
at $h=0$, see Prop. \ref{prop: subdif of HK}. 
\end{rmrk}
\section{Minimizing Movement approach}\label{sec: 3} 
\subsection{Theorem and Assumptions}\label{subsec: 3.1}
Let us return to the setting described in Sect. \ref{subsec: our MM approach} with $X:=\bar{\Omega}\subset\xR^d, \ {\sf d}(x_1, x_2):=|x_1-x_2|,$ and let us define $\Phi(\tau, \mu, \nu):= \mathcal{E}(\nu) + \frac{1}{2\tau}\HK_{\Lambda, \Sigma}(\nu, \mu)^2, \ \mathcal{E} := \mathcal{F} + \mathcal{V},$ as in \eqref{eq: our MM} and \eqref{eq: our E}. Our goal is to prove that the associated Minimizing Movement scheme yields weak solutions to the scalar reaction-diffusion equation \eqref{eq: our eq} with no-flux boundary condition \eqref{eq: our bc}. In order to find appropriate assumptions on $F: [0, +\infty)\to \xR$ and $V: \bar{\Omega}\to \xR$, it is worth taking a look at the natural coercivity assumptions which typically arise in connection with the Minimizing Movement approach to gradient flows (cf. the fundamental boook \cite{AmbrosioGigliSavare05} by Ambrosio, Gigli and Savar\'e). Note that, by Prop. \ref{prop: basics on HK}, $(\mathcal{M}(\bar{\Omega}), \HK_{\Lambda, \Sigma})$ is a complete metric space. 
\begin{prpstn}[Minimizing Movement approach to gradient flows \cite{AmbrosioGigliSavare05}]\label{prop: typical ass}
Let $(\mathscr{S}, d)$ be a complete metric space and apply the Minimizing Movement scheme \eqref{eq: MM scheme} to $\Phi(\tau, v, x):= \mathcal{E}(x) + \frac{1}{2\tau}d(x,v)^2$. We suppose that the functional $\mathcal{E}: \mathscr{S}\to (-\infty, +\infty]$ satisfies the following assumptions: 
\begin{enumerate}
\item [(A1)] There exist $A, B > 0, \ x_\star\in\mathscr{S}$ such that
\begin{equation}\label{eq: A1}
\mathcal{E}(\cdot) \ \geq \ -A-Bd(\cdot, x_\star)^2. 
\end{equation}
\item [(A2)] $\mathcal{E}$ is lower semicontinuous, i.e.
\begin{equation}\label{eq: A2}
d(x_n, x) \to 0 \quad \Rightarrow \quad \liminf_{n\to\infty}\mathcal{E}(x_n) \geq \mathcal{E}(x). 
\end{equation}
\item [(A3)] Every $d$-bounded set contained in a sublevel of $\mathcal{E}$ is relatively compact, i.e.
\begin{equation}\label{eq: A3}
\sup_{n, m} \{\mathcal{E}(x_n), d(x_n, x_m)\} < +\infty \quad \Rightarrow \quad \exists n_k\uparrow +\infty, x\in\mathscr{S}: \ d(x_{n_k}, x)\to 0. 
\end{equation}
\end{enumerate}
Then, for every $u_0\in\{\mathcal{E} < +\infty\}$, the set of Generalized Minimizing Movements $\mathrm{GMM}(\Phi; u_0)$ is nonempty. Moreover, every $u\in\mathrm{GMM}(\Phi; u_0)$ is continuous (locally absolutely continuous even) and satisfies the energy dissipation inequality \eqref{eq: energy inequality}. 
\end{prpstn}  
\begin{proof}
See Chaps. 1-3 in \cite{AmbrosioGigliSavare05}. The definitions asscociated with the energy dissipation inequality \eqref{eq: energy inequality} can be found therein, too. A brief outline of \eqref{eq: energy inequality} is given in Sect. \ref{subsec: 3.3} in this paper.  
\end{proof}
So let us break down the Assumptions (A1), (A2) and (A3) of Prop. \ref{prop: typical ass} on $\mathcal{E}$ into assumptions on $F$ and $V$. Let us put the focus on $F$ first, supposing that $V\equiv 0$. We start with Ass. (A1). We notice that for a Borel measurable function $F: [0, +\infty)\to \xR$, we have
\begin{eqnarray*}
\int_\Omega{\min\{F(u(x)), 0\}\xdif x} \ > \ -\infty 
\end{eqnarray*}
for all $u: \Omega\to[0, +\infty), \ u\in\xLone(\Omega),$ if and only if $F$ is linearly bounded from below, i.e. condition \eqref{eq: A1 F} below is a necessary and sufficient condition for the well-posedness of $\mathcal{F}$ 
(cf. Thm. 5.1 and Ex. 5.5 in \cite{fonseca2007modern}).  %and notice that, if \eqref{eq: A1} holds good for some $x_\star\in\mathscr{S}$, then it clearly holds good for any $x_\star\in\mathscr{S}$ (with modified constants $A, B >0$). %The next lemma says that $\mathcal{F}$ satisfies Ass. (A1) if and only if $F$ is linearly bounded from below. 
\begin{lmm}\label{lmm: A1}
We suppose that $F: [0, +\infty) \to \xR$ is Borel measurable and there exists $C_F > 0$ such that
\begin{equation}\label{eq: A1 F}
F(s) \ \geq \ -C_F s - C_F \quad \text{ for all } s\in [0, +\infty).
\end{equation} 
Let $\mathcal{F}$ be defined as in \eqref{eq: our E} and let $\eta_0$ denote the null measure. Then there exist $A, B >0$ such that
\begin{equation}\label{eq: A1 HK}
\mathcal{F}(\cdot) \ \geq \ -A-B\HK_{\Lambda, \Sigma}(\cdot, \eta_0)^2. 
\end{equation}
%On the other hand, if a functional $\mathcal{F}: \mathcal{M}(\bar{\Omega})\to (-\infty, +\infty]$ defined as in \eqref{eq: our E} satisfies \eqref{eq: A1 HK}, then \eqref{eq: A1 F} necessarily holds good for some $C_F > 0$. 
\end{lmm}    
\begin{proof}
Obviously, \eqref{eq: A1 F} implies \eqref{eq: A1 HK} with $A:= C_F\mathscr{L}^d(\Omega)$ and $B:=\frac{\Sigma}{4}C_F$; the only thing to note is that we have $\frac{\Sigma}{4}\HK_{\Lambda, \Sigma}(\mu, \eta_0)^2 = \int_{\Omega}{u(x)\xdif x}$ if $\mu=u\mathscr{L}^d$ by Prop. \ref{prop: basics on HK} and $\mathscr{L}^d(\Omega) < +\infty$ since $\Omega$ is bounded. %Let us prove the opposite direction by contradiction and suppose that $\mathcal{F}$ satisfies \eqref{eq: A1 HK} for some $A, B > 0$ but \eqref{eq: A1 F} does not hold good. Then for every $N\in\xN$ there exists $s_N\in[0, +\infty)$ such that
%\begin{equation*}
%F(s_N) \ < \ -Ns_N-N, 
%\end{equation*} 
%and for $\nu_N := s_N\mathscr{L}^d$, we have 
%\begin{equation*}
%\mathcal{F}(\nu_N) = F(s_N)\mathscr{L}^d(\Omega) < -Ns_N\mathscr{L}^d(\Omega)-N\mathscr{L}^d(\Omega) = -N\frac{\Sigma}{4}\HK_{\Lambda, \Sigma}(\nu_N, \eta_0)^2-N\mathscr{L}^d(\Omega),
%\end{equation*}
%which is a contradiction to \eqref{eq: A1 HK}. Lemma \ref{lmm: A1} is proved. 
\end{proof}
%Also condition \eqref{eq: A1 F} eventually justifies our definition of $\mathcal{F}$ as a mapping from $\mathcal{M}(\bar{\Omega})$ into $(-\infty, +\infty]$. Indeed,  

As $\HK_{\Lambda, \Sigma}$ metrizes the weak topology on $\mathcal{M}(\bar{\Omega})$ in duality with continuous and bounded functions, according to Prop. \ref{prop: basics on HK}, compactness issues such as \eqref{eq: A3} are closely linked to an application of Prokhorov's Theorem. Thus the compactness of $(\bar{\Omega}, |\cdot|)$ yields the relative compactness of every $\HK_{\Lambda, \Sigma}$-bounded set in $(\mathcal{M}(\bar{\Omega}), \HK_{\Lambda, \Sigma})$ (cf. Thm. 2.2, Cor. 7.16 in \cite{liero2018optimal}), i.e. \eqref{eq: A3} holds good in any case. However, this does not suffice in view of Ass. (A2) of Prop. \ref{prop: typical ass} and the fact that we aim to obtain (weak) solutions $u: [0, +\infty)\times \Omega\to [0, +\infty)$ to \eqref{eq: our eq}, \eqref{eq: our bc}. We need a condition on $F$ which guarantees
\begin{equation}\label{eq: A2 F}
\sup_n\{\mathcal{F}(\mu_n), \HK_{\Lambda, \Sigma}(\mu_n, \eta_0)\} < +\infty \quad \Rightarrow \quad \exists n_k\uparrow +\infty, \ u:\Omega\to [0, +\infty): \  \HK_{\Lambda, \Sigma}(\mu_{n_k}, u\mathscr{L}^d) \to 0.
\end{equation} 
    
\begin{lmm}\label{lmm: A2 F}
We suppose that $F$ is Borel measurable and linearly bounded from below \eqref{eq: A1 F}. Let $\mathcal{F}$ be defined as in \eqref{eq: our E}. Then \eqref{eq: A2 F} holds good if and only if $F$ has superlinear growth, i.e.
\begin{equation}\label{eq: superlinear growth}
\lim_{s\to\infty}\frac{F(s)}{s} \ = \ +\infty, 
\end{equation}
and it is equivalent to 
\begin{equation}\label{eq: 69}
\sup_n\{\mathcal{F}(\mu_n), \HK_{\Lambda, \Sigma}(\mu_n, \eta_0)\} < +\infty, \ \mu_n = u_n\mathscr{L}^d \quad \Rightarrow \quad \exists n_k\uparrow +\infty, \ u:\Omega\to [0, +\infty): \  u_{n_k}\stackrel{\xLone}{\rightharpoonup} u. 
\end{equation}
Furthermore, $\mathcal{F}$ is lower semicontinuous, i.e.
\begin{equation}
\HK_{\Lambda, \Sigma}(\mu_n, \mu) \to 0 \quad \Rightarrow \quad \liminf_{n\to\infty}\mathcal{F}(\mu_n)\geq\mathcal{F}(\mu),
\end{equation}
if and only if $F$ is convex, lower semicontinuous and has superlinear growth \eqref{eq: superlinear growth}. 
\end{lmm}
\begin{proof}
We suppose that $F$ has superlinear growth, $\sup_n\mathcal{F}(\mu_n) < +\infty, \ \mu_n=u_n\mathscr{L}^d$ and $\sup_n\int_{\Omega}{u_n(x)\xdif x} = \frac{\Sigma}{4}\sup_n\HK_{\Lambda, \Sigma}(\mu_n, \eta_0)^2 < +\infty$. By \eqref{eq: superlinear growth}, for every $M\in\xN$ there exists $s_M > 0$ such that $F(s) \geq Ms$ for all $s\geq s_M$. Let a Borel set $E\subset\Omega$ be given. We fix $n\in\xN$ and define $E_M:=E\cap\{x\in\Omega: \ u_n(x) > s_M\}, \ (M\in\xN)$. We have
\begin{eqnarray*}
\int_E{u_n(x)\xdif x} &\leq& s_M\mathscr{L}^d(E) + \frac{1}{M}\int_{E_M}{F(u_n(x))\xdif x} \\
&\leq& s_M\mathscr{L}^d(E) + \frac{1}{M}\Big(\mathcal{F}(u_n) + C_F\int_{\Omega}{u_n(x)\xdif x} + C_F\mathscr{L}^d(\Omega)\Big),
\end{eqnarray*} 
which shows the equiintegrability of $(u_n)_n$. Therefore, by Dunford-Pettis-Theorem, there exist a subsequence $n_k\uparrow +\infty$ and $u: \Omega\to [0, +\infty)$ such that $u_{n_k} \stackrel{\xLone}{\rightharpoonup} u$, where $\stackrel{\xLone}{\rightharpoonup}$ denotes weak convergence in $\xLone(\Omega)$. Clearly, $u_{n_k}\stackrel{\xLone}{\rightharpoonup}u$ implies $\HK_{\Lambda, \Sigma}(\mu_{n_k}, u\mathscr{L}^d)\to 0$.  

Now, let us suppose that $F$ does not have superlinear growth \eqref{eq: superlinear growth}. Then there exists $s_n\to\infty$ such that $\sup_n\frac{F(s_n)}{s_n} < +\infty$. We fix $\bar{x}\in\Omega$. Let $B_n$ be the open ball around $\bar{x}$ with radius $\Big(\frac{1}{s_n}\Big)^{1/d}$, and let $\mu_n := u_n\mathscr{L}^d, \ u_n(x):=s_n \text{ if } x\in\Omega\cap B_n, u_n(x):= 0 \text{ else }$. Then $\sup_n\{\mathcal{F}(\mu_n), \HK_{\Lambda, \Sigma}(\mu_n, \eta_0)\} < +\infty$ and $\mu_n$ converges to the Dirac measure $\mathrm{Vol}_d\delta_{\bar{x}}$ (where $\mathrm{Vol}_d$ denotes the volume of the unit ball in $\xR^d$) in $(\mathcal{M}(\bar{\Omega}), \HK_{\Lambda, \Sigma})$. Hence, \eqref{eq: A2 F} does not hold. 

So we have proved that \eqref{eq: A2 F} holds good if and only if $F$ has superlinear growth and it is equivalent to \eqref{eq: 69}. 

The second part of the lemma follows from the first part, \eqref{eq: A1 F} and the fact that, for a Borel measurable function $G: \xR\to [0, +\infty)$, the functional $\mathcal{G}: \xLone(\Omega)\to [0, +\infty], \ \mathcal{G}(u):= \int_{\Omega}{G(u(x))\xdif x},$  is lower semicontinuous w.r.t. weak $\xLone$-convergence if and only if $G$ is convex and lower semicontinuous (cf. Thms. 5.9 and 5.14 in \cite{fonseca2007modern}).  
\end{proof}

Note that, if $F$ is convex, then it is automatically linearly bounded from below, i.e. there exists $C_F > 0$ such that \eqref{eq: A1 F} holds good. We have seen so far that $\mathcal{F}$ satisfies the assumptions of Prop. \ref{prop: typical ass} if and only if $F$ is convex, lower semicontinuous and has superlinear growth \eqref{eq: superlinear growth}. In addition, the proof of (a weak form of) the reaction-diffusion equation \eqref{eq: our eq} with no-flux boundary condition \eqref{eq: our bc} will require a sort of differentiability property of $\mathcal{F}$, see Ass. \ref{ass: F} below. This condition will arise quite naturally. %Moreover, we will need a technical assumption which makes the passage to the limit $\tau_k\downarrow 0$ in the discrete weak version of \eqref{eq: our eq} possible, see condition \eqref{eq: technical ass} and Rem. \ref{rem: technical ass} below.  

Now, our theorem reads as follows.   

\begin{thrm}\label{thm: main thm}
Let a continuous strictly convex function $F: [0, +\infty)\to\xR$ with superlinear growth \eqref{eq: superlinear growth} and a Lipschitz continuous function $V: \bar{\Omega}\to\xR$ be given and define $\mathcal{E}:=\mathcal{F}+\mathcal{V}: \mathcal{M}(\bar{\Omega})\to (-\infty, +\infty]$ and $\Phi$ according to \eqref{eq: our E} and \eqref{eq: our MM}. Let $F$ be differentiable in $(0, +\infty)$ and define $L_F, \hat{L}_F: [0, +\infty) \to \xR$ as 
\begin{equation}\label{eq: LF}
L_F(s) := \begin{cases}
sF'(s) - F(s) &\text{ if } s\in (0, +\infty), \\
-F(0) &\text{ if } s=0,
\end{cases}
\quad \quad \hat{L}_F(s):=L_F(s) + F(s) =
\begin{cases}
sF'(s) &\text{ if } s\in (0, +\infty), \\
0 &\text{ if } s=0. 
\end{cases}
\end{equation}
We suppose that $\mathcal{F}$ satisfies Ass. \ref{ass: F} (see below). 

Then the following holds good: For every $\mu_0\in\{\mathcal{E} < +\infty\}$, the set $\mathrm{GMM}(\Phi; \mu_0)$ is nonempty. Furthermore, if $\mu\in\mathrm{GMM}(\Phi; \mu_0)$, then there exists a curve $u: [0, +\infty) \to \xLone(\Omega), \ u\geq 0,$ such that 
\begin{eqnarray}
&&\mu(t) \ = \ u(t)\mathscr{L}^d, \\
&&u(t_n) \ \stackrel{\xLone}{\rightharpoonup} \ u(t) \quad \text{ if } t_n\to t, 
\end{eqnarray}
for all $t\geq 0$, and $u$ is a solution to a weak version of the reaction-diffusion equation \eqref{eq: our eq} with no-flux boundary condition \eqref{eq: our bc}, i.e.  
\begin{equation}\label{eq: regularity} \hat{L}_F(u)\in\xLtwo_{\mathrm{loc}}([0,+\infty); \xLone(\Omega)), \quad\quad
%L_F(u(t))\in \rm{W}^{1,1}(\Omega) \quad \text{ for a.e. } t > 0, 
L_F(u)\in \xLtwo_{\mathrm{loc}}([0,+\infty); \rm{W}^{1,1}(\Omega))   
\end{equation}
and
\begin{equation}\label{eq: weak form 1}
\mathcal{I}_{F,V,\psi, u} \ = \ \int_0^\infty{\int_\Omega{u(t,x)\partial_t\psi(t,x)\xdif x}\xdif t} + \int_\Omega{u(0, x)\psi(0,x)\xdif x}
\end{equation} 
for all $\psi\in\xCtwo_{\mathrm{c}}(\xR\times\xR^d)$, where $\mathcal{I}_{F,V,\psi,u}$ is defined as
\begin{equation}\label{eq: weak form 2}
\int_0^\infty{\int_\Omega{\Big[\Lambda\Big\langle\nabla L_F(u(t, x)) + u(t,x)\nabla V(x), \nabla_x \psi(t,x)\Big\rangle + \Sigma \Big(\hat{L}_F(u(t,x)) +V(x)u(t,x)\Big)\psi(t,x)\Big]\xdif x}\xdif t}.
\end{equation}
\end{thrm}
\begin{rmrk}
$\xCtwo_{\mathrm{c}}(\xR\times \xR^d)$ denotes the set of all twice continuously differentiable functions $\psi: \xR\times\xR^d\to\xR$ with compact support in $\xR\times \xR^d$. In order to obtain a weak form of the scalar reaction-diffusion equation \eqref{eq: our eq}, it would suffice to prove \eqref{eq: weak form 1}, \eqref{eq: weak form 2} for all twice continuously differentiable functions $\psi: \xR\times \Omega\to\xR$ with compact support in $\xR\times\Omega$ (short $\psi\in\xCtwo_{\mathrm{c}}(\xR\times\Omega)$). Establishing \eqref{eq: weak form 1}, \eqref{eq: weak form 2} for all $\psi\in\xCtwo_{\mathrm{c}}(\xR\times\xR^d)$ instead means to include the no-flux boundary condition \eqref{eq: our bc} in a weak form, and will be an extra challenge in the proof of Thm. \ref{thm: main thm}. 

Moreover, $\nabla L_F(u(t, \cdot))\in\xLone(\Omega; \xR^d)$ denotes the weak gradient of $L_F(u(t))\in\rm{W}^{1,1}(\Omega)$, and $\nabla_x\psi$ denotes the gradient of $\psi$ with respect to the $x$-variable. %The well-posedness of the integral in \eqref{eq: weak form 2} will become apparent from the proof.  
\end{rmrk}
\begin{rmrk}\label{rem: technical ass}
By basic convex analysis, $L_F$ and $\hat{L}_F$ are continuous in $[0, +\infty)$ and $L_F$ is nondecreasing, i.e. $L_F(s_1)\leq L_F(s_2)$ whenever $s_1\leq s_2$. Furthermore, $L_F(s_n)\to +\infty$ if $s_n\to+\infty$ because $F$ has superlinear growth. Our assumption that $F$ is not only convex but strictly convex seems none too restrictive and makes things considerably easier since in this case, $L_F$ is strictly increasing, i.e. $L_F(s_1) < L_F(s_2)$ whenever $s_1 < s_2$. Thus, we have
\begin{equation}\label{eq: inverse LF}
L_F(s_n) \to L\in\xR \quad \quad \Rightarrow \quad \quad \exists s\in[0, +\infty): \ L=L_F(s), \quad s_n\to s.  
\end{equation}   
%Condition \eqref{eq: technical ass} will make the passage to the limit $\tau_k\downarrow 0$ in the reaction part possible. 
The fact that, according to \eqref{eq: inverse LF}, $\mathscr{L}^d$-a.e.-convergence of $L_F(u_n)$ %for $t>0$ and discrete solutions $\mu_{\tau_{k_l}}(t)=u_{\tau_{k_l}}(t)\mathscr{L}^d$ to our Minimizing Movement scheme 
for $u_n: \Omega\to[0, +\infty), \ u_n\in\xLone(\Omega) \ (n\in\xN),$
translates into $\mathscr{L}^d$-a.e. convergence of $u_n$ will be a useful ingredient in our proof, cf. Sect. \ref{subsec: 3.2}.   
\end{rmrk}

Now, let us be precise about the differentiability condition imposed on $\mathcal{F}$.
\begin{ass}\label{ass: F}
We suppose that $L_F(u_0)\in\xLone(\Omega)$ and 
\begin{equation}\label{eq: diff cond}
\lim_{h\to 0} \frac{\mathcal{F}(\nu_h)-\mathcal{F}(\nu_0)}{h} \ = \ \int_{\Omega}{\Big[-L_F(u_0(x))\xtr \xDif v(x) + 2\hat{L}_F(u_0(x))R(x)\Big]\xdif x}, 
\end{equation}
whenever $v: \Omega\to \xR^d$ is continuously differentiable and has compact support in $\Omega$, $R: \Omega\to\xR$ is bounded and Borel measurable, $\nu_0 = u_0\mathscr{L}^d\in\{\mathcal{F} < +\infty\}$ and the curve $\mathcal{N}\ni h \mapsto \nu_h\in\mathcal{M}(\bar{\Omega})$ is defined according to \eqref{eq: our direction}, i.e. 
\begin{equation}
\nu_h:= (I+hv)_{\#}(1+hR)^2\nu_0 
\end{equation}
(where $L_F, \hat{L}_F$ are defined as in \eqref{eq: LF}, $\xDif v$ denotes the differential of $v$ and $\xtr \xDif v$ its trace). 
\end{ass}
A similar condition has already been treated in the study of diffusion equations \eqref{eq: diffusion equation} (cf. Sect. 10.4.3 in \cite{AmbrosioGigliSavare05}). The differentiability of $h\mapsto\mathcal{F}(\nu_h)$ at $h=0$, for such curves $h\mapsto\nu_h$, together with our analysis from Sect. \ref{subsec: directional derivative of HK}, will form the very basis for proving \eqref{eq: weak form 1}, \eqref{eq: weak form 2}. We note that, if $\nu_0=u_0\mathscr{L}^d$ and $\nu_h$ is defined as above, then,  
\begin{equation}\label{eq: 78}
\nu_h = u_h\mathscr{L}^d, \quad\quad \mathop{\rm det\,}\nolimits(\mathbb{I} + h\xDif v(x))u_h(x+hv(x)) \ = \ (1+hR(x))^2 u_0(x), \quad x\in\Omega, 
\end{equation} 
for $h$ in a neighbourhood around $0$ (where $\mathbb{I}$ denotes the identity matrix). This follows from the change of variables formula and the fact that, for $|h|$ small enough, $I+hv$ is a diffeomorphism mapping $\Omega$ onto $\Omega$ with $\mathop{\rm det\,}\nolimits(\xDif (I+hv)) = \mathop{\rm det\,}\nolimits(\mathbb{I} + h\xDif v) > 0$.
Moreover, for every $x\in\Omega$, the mapping $h\mapsto \mathop{\rm det\,}\nolimits(\mathbb{I} + h\xDif v(x))$ is differentiable at $h=0$ with derivative equal to $\xtr \xDif v(x)$. By \eqref{eq: 78} and the change of variables formula, we have 
\begin{equation}\label{eq: left-hand side}
\frac{\mathcal{F}(\nu_h) - \mathcal{F}(\nu_0)}{h} \ = \ \int_{\Omega}{\frac{1}{h}\Big[F\Big(\frac{(1+hR(x))^2u_0(x)}{\mathop{\rm det\,}\nolimits(\mathbb{I} + h\xDif v(x))}\Big)\mathop{\rm det\,}\nolimits(\mathbb{I} + h\xDif v(x)) - F(u_0(x))\Big]\xdif x}
\end{equation}  
if $|h|$ is small and $\nu_h\in\{\mathcal{F} < +\infty\}$. It is not difficult to see that the integrands of \eqref{eq: left-hand side} converge pointwise to the integrand of the right-hand side of \eqref{eq: diff cond} as $h\to 0$. So if the corresponding integrals also converge (e.g. by dominated convergence theorem or monotone convergence theorem), then \eqref{eq: diff cond} holds good.    
\begin{xmpl}\label{ex: 3.8}
We give two examples of functions $F: [0, +\infty)\to\xR$ satisfying the assumptions of Thm. \ref{thm: main thm}.
 
The first example is \begin{equation*}
F(s) := \begin{cases}
c_1s\log s &\text{ if } s\in (0, +\infty), \\
0 &\text{ if } s=0, 
\end{cases}
\quad \quad (c_1 > 0),
\end{equation*}
for which Thm. \ref{thm: main thm} yields \eqref{eq: weak form 1}, \eqref{eq: weak form 2} with
\begin{equation*}
L_F(s) = c_1s, \quad \hat{L}_F(s) = \begin{cases} c_1s + c_1s\log s &\text{ if } s\in (0, +\infty), \\ 0 &\text{ if } s=0. \end{cases}
\end{equation*}
In this case, Ass. \ref{ass: F} is established by simplifying the right-hand side of \eqref{eq: left-hand side} to 
\begin{eqnarray*} 
\int_{\Omega}{\Big[F(u_0(x))\frac{(1+hR(x))^2-1}{h} + c_1 u_0(x)(1+hR(x))^2 \ \frac{\log (1+hR(x))^2 - \log \mathop{\rm det\,}\nolimits(\mathbb{I} + h\xDif v(x))}{h}\Big]\xdif x} 
%+ \int_\Omega{c_2 u_0(x)^p\frac{(1+hR(x))^{2p}(\mathop{\rm det\,}\nolimits(\mathbb{I} + h\xDif v(x)))^{1-p} -1}{h}\xdif x}
\end{eqnarray*}
(using arithmetical rules of the logarithm) and by applying the dominated convergence theorem (note that $\int_{\Omega}{|F(u_0(x))|\xdif x} < +\infty$ if $\nu_0 = u_0\mathscr{L}^d\in\{\mathcal{F} < +\infty\}$). %Moreover, it is not difficult to see that $|s\log s|\leq \frac{1}{e} + s^p$ so that \eqref{eq: technical ass} holds good for $C:=\frac{c_1+c_2}{c_2(p-1)} + \frac{c_1}{e}$.  

The second basic example to which Thm. \ref{thm: main thm} is applicable is 
\begin{equation*}
F(s) := -c_1s^q + c_2 s^p \quad (c_1 \geq 0, \ c_2 > 0, \ p>1, \ q\in (0,1))
\end{equation*}
with 
\begin{equation*}
L_F(s)=c_1(1-q)s^q+c_2(p-1)s^p, \quad \hat{L}_F(s)=-c_1q s^q+c_2ps^p. 
\end{equation*}
In this case, the right-hand side of \eqref{eq: left-hand side} reads as
\begin{equation*}
\int_\Omega{\Big[-c_1 u_0(x)^q \ \frac{(\mathop{\rm det\,}\nolimits(\mathbb{I} + h\xDif v(x)))^{1-q} (1+hR(x))^{2q} - 1}{h} + c_2 u_0(x)^p\frac{(1+hR(x))^{2p}(\mathop{\rm det\,}\nolimits(\mathbb{I} + h\xDif v(x)))^{1-p} -1}{h}\Big]\xdif x},
\end{equation*}
and again, Ass. \ref{ass: F} can be established by using the dominated convergence theorem (note that $\int_\Omega{u_0(x)^q\xdif x} \leq \mathscr{L}^d(\Omega) + \int_\Omega{u_0(x)\xdif x}$). %Moreover, $F$ satisfies condition \eqref{eq: technical ass} with $C:=\frac{1}{\min\{1-q,p-1\}}$. 
\end{xmpl}
Finally, we remark that the Lipschitz continuity of $V$ appears to be a convenient condition for our purposes since in this case, $V\in\xCzero_b(\bar{\Omega})$ and thus $\mathcal{V}$ (defined as in \eqref{eq: our E}) obviously satisfies Ass. (A1), (A2), (A3) of Prop. \ref{prop: typical ass}, and in addition, $V$ is a.e. differentiable with bounded gradient $\nabla V$, and, by dominated convergence theorem,  
\begin{equation}\label{eq: diff cond V}
\lim_{h\to 0}\frac{\mathcal{V}(\nu_h)-\mathcal{V}(\nu_0)}{h} \ = \ \int_{\Omega}{\Big[\langle\nabla V(x), v(x)\rangle + 2V(x)R(x)\Big]u_0(x)\xdif x}
\end{equation}
for every curve $\mathcal{N}\ni h \mapsto \nu_h:=(I+hv)_{\#}(1+hR)^2\nu_0\in\mathcal{M}(\bar{\Omega}), \ \nu_0 = u_0\mathscr{L}^d,$ defined according to Ass. \ref{ass: F}. 
\subsection{Proof}\label{subsec: 3.2}
We prove Thm. \ref{thm: main thm}. 
\begin{proof}
If the assumptions of Thm. \ref{thm: main thm} hold, then Lem. \ref{lmm: A1}, Lem. \ref{lmm: A2 F} and the discussion in Sect. \ref{subsec: 3.1} show that Prop. \ref{prop: typical ass} is applicable to $\Phi$ and $\mathcal{E}:=\mathcal{F}+\mathcal{V}: \mathcal{M}(\bar{\Omega})\to (-\infty, +\infty]$ defined as in \eqref{eq: our MM} and \eqref{eq: our E}. Hence, for every $\mu_0=u_0\mathscr{L}^d\in\{\mathcal{E} < +\infty\}$, the set $\mathrm{GMM}(\Phi; \mu_0)$ is nonempty. So let $\mu\in\mathrm{GMM}(\Phi; \mu_0)$. There exist a subsequence $(\tau_k)_{k\in\xN}, \tau_k\downarrow 0,$ and discrete solutions $\mu_{\tau_k}$ to the associated Minimizing Movement scheme \eqref{eq: MM scheme} converging pointwise to $\mu$, i.e. $\mu_{\tau_k}(0) \ = \ \mu_0$ and 
\begin{equation*}
\lim_{k\to\infty}\HK_{\Lambda, \Sigma}(\mu_{\tau_k}(t), \mu(t)) \ = \ 0 \quad \quad \text{ for all } t\geq 0.  
\end{equation*}
Every discrete solution $\mu_{\tau_k}$ is assigned a curve $u_{\tau_k}: [0, +\infty)\to \xLone(\Omega)$ such that $u_{\tau_k}(0)=u_0, \ u_{\tau_k}(t, \cdot)\geq 0$ and 
\begin{equation*}
\mu_{\tau_k}(t) \ = \ u_{\tau_k}(t)\mathscr{L}^d \quad \quad \text{ for all } t\geq 0. 
\end{equation*} 
We note that, by \eqref{eq: MM scheme}, $t\mapsto\mathcal{E}(\mu_{\tau_k}(t))$ is decreasing. Since $V\in\xCzero_b(\bar{\Omega})$ and $\sup_k \HK_{\Lambda, \Sigma}(\mu_{\tau_k}(t), \eta_0) < +\infty$, it follows from $\mathcal{E}(\mu_{\tau_k}(t))\leq \mathcal{E}(\mu_0)$ that $\sup_{k}\mathcal{F}(\mu_{\tau_k}(t)) < +\infty$. Thus, according to Lem. \ref{lmm: A2 F}, there exists a curve $u: [0, +\infty) \to \xLone(\Omega)$ such that $u(0) = u_0, \ u(t, \cdot)\geq 0$ and 
\begin{equation*}
u_{\tau_k}(t) \ \stackrel{\xLone}{\rightharpoonup} \ u(t), \quad \quad \mu(t) \ = \ u(t)\mathscr{L}^d \quad \quad \text{ for all } t\geq0
\end{equation*} 
(this convergence holds good for the whole sequence $\tau_k\downarrow 0$ as we already know that $\HK_{\Lambda, \Sigma}(\mu_{\tau_k}(t), \mu(t))\to 0$). By Prop. \ref{prop: typical ass}, $\mu$ is continuous, i.e. $\HK_{\Lambda, \Sigma}(\mu(t_n), \mu(t))\to 0$ whenever $t_n\to t$. Since $\mathcal{E}(\mu(t))\leq \mathcal{E}(\mu_0)$ for all $t\geq 0$, the same arguments as before show that
\begin{equation*}
u(t_n) \ \stackrel{\xLone}{\rightharpoonup} \ u(t) \quad \quad \text{ whenever } t_n\to t, \ t\geq 0.  
\end{equation*} 

Now, let $v: \Omega\to\xR^d$ be a continuously differentiable function with compact support in $\Omega$ and let $R: \Omega\to\xR$ be a bounded Borel measurable function. We define, for $k\in\xN, \ n\in\xN$, the curve
\begin{equation*}
\mathcal{N}\ni h \mapsto \nu_h:=(I+hv)_{\#}(1+hR)^2\mu_{\tau_k}(n\tau_k)\in\mathcal{M}(\bar{\Omega})
\end{equation*}
according to \eqref{eq: our direction}. We recall that
\begin{equation*}
\mu_{\tau_k}(n\tau_k) \text{ is a minimizer for } \Phi(\tau_k, \mu_{\tau_k}((n-1)\tau_k), \cdot) = \mathcal{E}(\cdot) + \frac{1}{2\tau_k}\HK_{\Lambda, \Sigma}(\cdot, \mu_{\tau_k}((n-1)\tau_k))^2,
\end{equation*} 
and we establish a necessary condition of first order involving the Fr\'echet subdifferential of 
\begin{equation*}
h\mapsto -\frac{1}{2}\HK_{\Lambda, \Sigma}(\nu_h, \mu_{\tau_k}((n-1)\tau_k))^2
\end{equation*}
at $h=0$ and the directional derivatives \eqref{eq: diff cond} and \eqref{eq: diff cond V} of $\mathcal{F}$ and $\mathcal{V}$. We set $\mu_{\tau_k}^n:=\mu_{\tau_k}(n\tau_k), \ u_{\tau_k}^n:=u_{\tau_k}(n\tau_k)$. Let $\beta_{\tau_k}^n\in\mathcal{M}(\mathfrak{C}\times\mathfrak{C})$ be optimal in the definition of $\HK_{\Lambda, \Sigma}(\mu_{\tau_k}^n, \mu_{\tau_k}^{n-1})^2$ according to \eqref{eq: cone 1}, \eqref{eq: Wasserstein on the cone}, with first marginal $\alpha_{\tau_k}^n\in\mathcal{M}_2(\mathfrak{C})$, $\mathfrak{h}\alpha_{\tau_k}^n = \mu_{\tau_k}^n$, and second marginal $\alpha_{\tau_k}^{n-1}, \ \mathfrak{h}\alpha_{\tau_k}^{n-1} = \mu_{\tau_k}^{n-1}$. Since $\mu_{\tau_k}^n$ is a minimizer for $\Phi(\tau_k, \mu_{\tau_k}^{n-1}, \cdot)$, we have
\begin{equation*}
\frac{\mathcal{E}(\nu_h) - \mathcal{E}(\mu_{\tau_k}^n)}{h} \ \geq \ \frac{\frac{1}{2\tau_k}\HK_{\Lambda, \Sigma}(\mu_{\tau_k}^n, \mu_{\tau_k}^{n-1})^2 - \frac{1}{2\tau_k}\HK_{\Lambda, \Sigma}(\nu_h, \mu_{\tau_k}^{n-1})^2}{h},  \quad \quad h\in\mathcal{N}, \ h > 0,  
\end{equation*} 
and passing to the limit $h\downarrow 0$, we obtain 
\begin{eqnarray*}
\int_{\Omega}{\Big[-L_F(u_{\tau_k}^n(x))\xtr \xDif v(x) + 2\hat{L}_F(u_{\tau_k}^n(x))R(x) + \Big(\langle\nabla V(x), v(x)\rangle + 2V(x)R(x)\Big)u_{\tau_k}^n(x)\Big]\xdif x} \ \geq \ \frac{1}{\tau_k} \ \mathfrak{F}_{\tau_k, n, v, R}, 
\end{eqnarray*}
with
\begin{eqnarray*}
\mathfrak{F}_{\tau_k, n, v, R}:= \frac{4}{\Sigma}\int_{\mathfrak{C}\times\mathfrak{C}}{\Big[-r_1^2R(x_1)+r_1r_2R(x_1)\cos(\sqrt{\Sigma/4\Lambda} \ |x_1-x_2|)+r_1r_2\sqrt{\Sigma/4\Lambda} \ \langle S_{\Lambda,\Sigma}(x_1,x_2), v(x_1)\rangle\Big]\xdif\beta_{\tau_k}^n}, \\
\quad\quad\quad S_{\Lambda, \Sigma}(x_1, x_2):= \begin{cases}
\frac{\sin(\sqrt{\Sigma/4\Lambda} \ |x_1-x_2|)}{|x_1-x_2|}(x_2-x_1) &\text{ if } x_1\neq x_2, \\
0 &\text{ if } x_1=x_2, 
\end{cases}
\end{eqnarray*}
by Prop. \ref{prop: subdif of HK}, Ass. \ref{ass: F} and \eqref{eq: diff cond V}. As we can switch between $v, R$ and $-v, -R$, the following necessary condition of first order holds good
\begin{equation}\label{eq: cond of first order}
\int_{\Omega}{\Big[-L_F(u_{\tau_k}^n)\xtr \xDif v + 2\hat{L}_F(u_{\tau_k}^n)R + \Big(\langle\nabla V, v\rangle + 2VR\Big)u_{\tau_k}^n\Big]\xdif x} \ = \ \frac{1}{\tau_k} \ \mathfrak{F}_{\tau_k, n, v, R}
\end{equation}  
for all continuously differentiable functions $v: \Omega\to\xR^d$ with compact support in $\Omega$ and all bounded Borel measurable functions $R: \Omega\to\xR$.

Let $\phi: \Omega\to \xR$ be a twice continuously differentiable function with compact support in $\Omega$. Setting $v:=\frac{4\Lambda}{\Sigma}\nabla\phi$, $R:=2\phi$ and applying our necessary condition of first order \eqref{eq: cond of first order} and Prop. \ref{prop: connection of subdif with MM}, we obtain
\begin{eqnarray*}
&& \Big|\int_\Omega{\Big[-\Lambda L_F(u_{\tau_k}^n)\Delta\phi + \Sigma\hat{L}_F(u_{\tau_k}^n)\phi + \Big(\Lambda\langle \nabla V, \nabla\phi\rangle + \Sigma V \phi\Big)u_{\tau_k}^n\Big]\xdif x} - \frac{1}{\tau_k}\Big(\int_\Omega{\phi\xdif\mu_{\tau_k}^{n-1}} - \int_{\Omega}{\phi\xdif \mu_{\tau_k}^n}\Big)\Big| \\ &\leq& \frac{1}{\tau_k} C_{\phi,\Lambda, \Sigma} \HK_{\Lambda, \Sigma}(\mu_{\tau_k}^n, \mu_{\tau_k}^{n-1})^2, 
\end{eqnarray*} 
with $C_{\phi, \Lambda, \Sigma} \geq 0$ only depending on $\phi, \Lambda, \Sigma$ (cf. \eqref{eq: 54}, \eqref{eq: 55}). Hence, for every $\psi\in\xCtwo_{\mathrm{c}}(\xR\times \Omega)$ there exist $C_{\psi, \Lambda, \Sigma} > 0, \ N_{\psi, \tau_k}\in\xN$ such that $\psi(t, \cdot)\equiv0$ for $t\geq N_{\psi, \tau_k}\tau_k, \ \sup_k N_{\psi, \tau_k}\tau_k < +\infty$ and
\begin{eqnarray*}
&&\Big|\int_0^\infty{\mathfrak{I}_{\tau_k, \psi}(t)\xdif t} - \int_{\tau_k}^{\infty}{\int_\Omega{\frac{u_{\tau_k}(t-\tau_k, x) - u_{\tau_k}(t, x)}{\tau_k}\psi(t,x)\xdif x}\xdif t} - \int_{0}^{\tau_k}{\int_\Omega{\frac{u_{\tau_k}(0, x) - u_{\tau_k}(t, x)}{\tau_k}\psi(t,x)\xdif x}\xdif t}\Big| \\ && \leq \ C_{\psi, \Lambda, \Sigma}\sum_{n=1}^{N_{\psi, \tau_k}}{\HK_{\Lambda, \Sigma}(\mu_{\tau_k}^n, \mu_{\tau_k}^{n-1})^2},
\end{eqnarray*}
where $\mathfrak{I}_{\tau_k, \psi}(t)$ is defined as 
\begin{equation*}
\int_\Omega{\Big[-\Lambda L_F(u_{\tau_k}(t,x))\Delta_x\psi(t,x) + \Sigma\hat{L}_F(u_{\tau_k}(t,x))\psi(t,x) + \Big(\Lambda\langle \nabla V(x), \nabla_x\psi(t,x)\rangle + \Sigma V(x) \psi(t,x)\Big)u_{\tau_k}(t,x)\Big]\xdif x}.
\end{equation*}
Standard tools from the theory of the Minimizing Movement approach to gradient flows yield 
\begin{eqnarray}
&&\sup \Big\{\HK_{\Lambda, \Sigma}(\mu_{\tau_k}^n, \eta_0)^2: \ 1\leq n \leq N_k, \ k\in\xN\Big\} < +\infty, \label{eq: standard tools 1} \\ &&\sum_{n=1}^{N_k}{\HK_{\Lambda, \Sigma}(\mu_{\tau_k}^n, \mu_{\tau_k}^{n-1})^2} \ \leq \ 2\tau_k\Big(\mathcal{E}(\mu_0) + A + B\HK_{\Lambda, \Sigma}(\mu_{\tau_k}^{N_k}, \eta_0)^2\Big) \quad \to 0 \quad\quad \text{ as } k\to\infty, \label{eq: standard tools 2} 
\end{eqnarray}
whenever $\sup_k N_k\tau_k < +\infty$ (and where $A, B > 0$ s.t. $\mathcal{E}(\cdot)\geq -A-B\HK_{\Lambda, \Sigma}(\cdot, \eta_0)^2$), see e.g. the first part of the proof in Sect. 3.2 in \cite{fleissner2016gamma}. Furthermore, we have
\begin{eqnarray*}
&&\int_{\tau_k}^{\infty}{\int_\Omega{\frac{u_{\tau_k}(t-\tau_k, x) - u_{\tau_k}(t, x)}{\tau_k}\psi(t,x)\xdif x}\xdif t} + \int_{0}^{\tau_k}{\int_\Omega{\frac{u_{\tau_k}(0, x) - u_{\tau_k}(t, x)}{\tau_k}\psi(t,x)\xdif x}\xdif t} \\
&& = \ \int_0^\infty{\int_\Omega{\frac{\psi(t+\tau_k, x) - \psi(t,x)}{\tau_k}u_{\tau_k}(t,x)\xdif x}\xdif t} + \frac{1}{\tau_k}\int_0^{\tau_k}{\int_\Omega{u_{\tau_k}(0,x)\psi(t,x)\xdif x}\xdif t} \\
&& \to \ \int_0^\infty{\int_\Omega{u(t,x)\partial_t\psi(t,x)\xdif x}\xdif t} + \int_\Omega{u(0,x)\psi(0,x)\xdif x} \quad \quad \text{ as } \tau_k\downarrow 0
\end{eqnarray*}
(since $\partial_t\psi$ is uniformly continuous and bounded, $u_{\tau_k}(t)\stackrel{\xLone}{\rightharpoonup} u(t)$, $\sup\Big\{\int_\Omega{u_{\tau_k}(t,x)\xdif x}: \ t\leq N_{\psi, \tau_k}\tau_k, \ k\in\xN\Big\} = \sup \Big\{\frac{\Sigma}{4}\HK_{\Lambda, \Sigma}(\mu_{\tau_k}(t), \eta_0)^2: \ t\leq N_{\psi, \tau_k}\tau_k, \ k\in\xN\Big\} < +\infty$). So in order to establish \eqref{eq: weak form 1} for $\psi\in\xCtwo_{\mathrm{c}}(\xR\times\Omega)$, all that remains is to prove that $L_F(u(t))\in\rm{W}^{1,1}(\Omega)$ for a.e. $t > 0$ and 
$\int_0^\infty{\mathfrak{I}_{\tau_k, \psi}(t)\xdif t} \to \mathcal{I}_{F,V,\psi, u}$ as $\tau_k\downarrow 0$. Again, the necessary condition of first order will smooth the way. We set $R\equiv 0$ in \eqref{eq: cond of first order} and obtain 
\begin{eqnarray*}
&&\Big|\int_\Omega{-L_F(u_{\tau_k}^n(x))\xtr \xDif v(x)\xdif x}\Big| \\ &&\leq \sup_{x\in\Omega}|\nabla V(x)| \int_{\Omega}{|v(x)|u_{\tau_k}^n(x)\xdif x} +\frac{1}{\tau_k} \Big(\int_{\mathfrak{C}\times\mathfrak{C}}{\frac{4}{\Sigma\Lambda}r_2^2|S_{\Lambda, \Sigma}(x_1,x_2)|^2\xdif\beta_{\tau_k}^n}\Big)^{1/2} \Big(\int_{\mathfrak{C}\times\mathfrak{C}}{r_1^2|v(x_1)|^2\xdif \beta_{\tau_k}^n}\Big)^{1/2} \\
&&\leq \Big(\sup_{x\in\Omega}|\nabla V(x)| \Big(\int_\Omega{u_{\tau_k}^n(x)\xdif x}\Big)^{1/2} + \frac{1}{\sqrt{\Lambda}}\frac{\HK_{\Lambda, \Sigma}(\mu_{\tau_k}^n, \mu_{\tau_k}^{n-1})}{\tau_k}\Big)\Big(\int_\Omega{|v(x)|^2u_{\tau_k}^n(x)\xdif x}\Big)^{1/2} \\
&&\leq  \Big(\sup_{x\in\Omega}|\nabla V(x)| \Big(\int_\Omega{u_{\tau_k}^n(x)\xdif x}\Big)^{1/2} + \frac{1}{\sqrt{\Lambda}}\frac{\HK_{\Lambda, \Sigma}(\mu_{\tau_k}^n, \mu_{\tau_k}^{n-1})}{\tau_k}\Big)\Big(\int_\Omega{u_{\tau_k}^n(x)\xdif x}\Big)^{1/2} \ \sup_{x\in\Omega}|v(x)| 
\end{eqnarray*} 
by applying H\"older's inequality and \eqref{eq: first derivatives R and theta} and \eqref{eq: right derivative at 0} which yield 
\begin{eqnarray*}
\int_{\mathfrak{C}\times\mathfrak{C}}{\frac{4}{\Sigma}r_2^2|S_{\Lambda, \Sigma}(x_1,x_2)|^2\xdif\beta_{\tau_k}^n} \ \leq \ \int_{\mathfrak{C}\times\mathfrak{C}}{{\sf d}_{\mathfrak{C},\Lambda, \Sigma}([x_1,r_1],[x_2,r_2])^2\xdif\beta_{\tau_k}^n} \ = \ \HK_{\Lambda, \Sigma}(\mu_{\tau_k}^n, \mu_{\tau_k}^{n-1})^2
\end{eqnarray*}
(note that \eqref{eq: first derivatives R and theta} holds good for $t=0$, too, as $\theta'_+(0) = \lim_{t\downarrow 0}\theta'(t)$ and $\mathcal{R}'_+(0) = \lim_{t\downarrow 0}\mathcal{R}'(t)$). Firstly, the above estimations (which hold for every continuously differentiable function $v: \Omega\to\xR^d$ with compact support in $\Omega$) show that $L_F(u_{\tau_k}^n)$ is a function of bounded variation in $\Omega$ (recall that $L_F(u_{\tau_k}^n)\in\xLone(\Omega)$ by Ass. \ref{ass: F}), and in a second step, they show, according to Hahn-Banach Theorem and the fact that the dual space of $\xLtwo(\Omega; u_{\tau_k}^n\mathscr{L}^d)$ can be identified with $\xLtwo(\Omega; u_{\tau_k}^n\mathscr{L}^d)$ itself, the existence of a function $w_{\tau_k}^n: \Omega \to \xR^d$ such that
\begin{eqnarray*}
\int_\Omega{|w_{\tau_k}^n(x)|^2u_{\tau_k}^n(x)\xdif x} &\leq& \Big(\sup_{x\in\Omega}|\nabla V(x)| \Big(\int_\Omega{u_{\tau_k}^n(x)\xdif x}\Big)^{1/2} + \frac{1}{\sqrt{\Lambda}}\frac{\HK_{\Lambda, \Sigma}(\mu_{\tau_k}^n, \mu_{\tau_k}^{n-1})}{\tau_k}\Big)^2 \\
\int_\Omega{-L_F(u_{\tau_k}^n(x))\xtr \xDif v(x)\xdif x} &=& \int_\Omega{\langle w_{\tau_k}^n(x), v(x)\rangle u_{\tau_k}^n(x)\xdif x}
\end{eqnarray*}    
for every continuously differentiable function $v: \Omega \to \xR^d$ with compact support in $\Omega$. This means that $L_F(u_{\tau_k}^n)\in\rm{W}^{1,1}(\Omega)$ with weak gradient $\nabla L_F(u_{\tau_k}^n) = w_{\tau_k}^n u_{\tau_k}^n$ and 
\begin{equation*}
\int_\Omega{\Big|\nabla L_F(u_{\tau_k}^n(x))\Big|\xdif x} \ \leq \ \Big(\sup_{x\in\Omega}|\nabla V(x)| \Big(\int_\Omega{u_{\tau_k}^n(x)\xdif x}\Big)^{1/2} + \frac{1}{\sqrt{\Lambda}}\frac{\HK_{\Lambda, \Sigma}(\mu_{\tau_k}^n, \mu_{\tau_k}^{n-1})}{\tau_k}\Big) \Big(\int_\Omega{u_{\tau_k}^n(x)\xdif x}\Big)^{1/2}.
\end{equation*}
Now, let us give an upper bound for $L_F(u_{\tau_k}^n)$ in $\xLone(\Omega)$ so that in the end, we have an upper bound for $L_F(u_{\tau_k}^n)$ in $\rm{W}^{1,1}(\Omega)$. We note that
\begin{equation*}
\int_\Omega{|L_F(u_{\tau_k}^n)|\xdif x} \ \leq \ \int_\Omega{\Big(|\hat{L}_F(u_{\tau_k}^n)| + |F(u_{\tau_k}^n)|\Big)\xdif x} \ \leq \ \int_\Omega{|\hat{L}_F(u_{\tau_k}^n)|\xdif x} + \mathcal{F}(\mu_{\tau_k}^n) + 2\int_\Omega{[|F(1)-F'(1)| + |F'(1)|u_{\tau_k}^n]\xdif x},
\end{equation*} 
and setting $v\equiv 0$ in \eqref{eq: cond of first order}, we obtain
\begin{eqnarray*}
\Big|\int_\Omega{\hat{L}_F(u_{\tau_k}^n(x))R(x)\xdif x}\Big| &\leq& \Big(\sup_{x\in\Omega}|V(x)|\int_\Omega{u_{\tau_k}^n(x)\xdif x}\Big) \sup_{x\in\Omega}|R(x)| \\ && + \ \frac{2}{\tau_k\Sigma}\Big(\int_{\mathfrak{C}\times\mathfrak{C}}{Q_{\Lambda, \Sigma}(x_1,r_1,x_2,r_2)^2\xdif \beta_{\tau_k}^n}\Big)^{1/2} \Big(\int_\Omega{u_{\tau_k}^n(x)\xdif x}\Big)^{1/2}\sup_{x\in\Omega}|R(x)|
\end{eqnarray*}
for all bounded Borel functions $R: \Omega\to\xR$, with $Q_{\Lambda, \Sigma}([x_1,r_1],[x_2,r_2]):=-r_1 + r_2\cos(\sqrt{\Sigma/4\Lambda}|x_1-x_2|)$ and 
\begin{equation*}
\int_{\mathfrak{C}\times\mathfrak{C}}{Q_{\Lambda, \Sigma}([x_1,r_1],[x_2,r_2])^2\xdif \beta_{\tau_k}^n} \ \leq \ \frac{\Sigma}{4}\int_{\mathfrak{C}\times\mathfrak{C}}{{\sf d}_{\mathfrak{C}, \Lambda, \Sigma}([x_1,r_1],[x_2,r_2])^2\xdif\beta_{\tau_k}^n} \ = \ \frac{\Sigma}{4}\HK_{\Lambda, \Sigma}(\mu_{\tau_k}^n, \mu_{\tau_k}^{n-1})^2
\end{equation*}
by \eqref{eq: first derivatives R and theta} (which also holds good at $t=0$) and \eqref{eq: right derivative at 0}. Hence, 
\begin{equation*}
\int_\Omega{|\hat{L}_F(u_{\tau_k}^n(x))|\xdif x} \ \leq \ \sup_{x\in\Omega}|V(x)|\int_\Omega{u_{\tau_k}^n(x)\xdif x} + \frac{1}{\sqrt{\Sigma}}\frac{\HK_{\Lambda, \Sigma}(\mu_{\tau_k}^n, \mu_{\tau_k}^{n-1})}{\tau_k} \Big(\int_\Omega{u_{\tau_k}^n(x)\xdif x}\Big)^{1/2}.
\end{equation*} 
Similarly, we obtain
\begin{equation*}
\int_E{|\hat{L}_F(u_{\tau_k}^n(x))|\xdif x} \ \leq \ \sup_{x\in\Omega}|V(x)|\int_E{u_{\tau_k}^n(x)\xdif x} + \frac{1}{\sqrt{\Sigma}}\frac{\HK_{\Lambda, \Sigma}(\mu_{\tau_k}^n, \mu_{\tau_k}^{n-1})}{\tau_k} \Big(\int_E{u_{\tau_k}^n(x)\xdif x}\Big)^{1/2}
\end{equation*} 
for every Borel set $E\subset \Omega$. 
Define $|\mu_{\tau_k}'|: (0, +\infty) \to [0, +\infty)$ as 
\begin{equation*}
|\mu_{\tau_k}'|(t):= \frac{\HK_{\Lambda, \Sigma}(\mu_{\tau_k}^n, \mu_{\tau_k}^{n-1})}{\tau_k} \quad \text{ for } t\in ((n-1)\tau_k, n\tau_k] \quad (n\in\xN). 
\end{equation*}
We have found out so far that $L_F(u_{\tau_k}(t))\in\rm{W}^{1,1}(\Omega)$ for all $t > 0$ and $k\in\xN$ and that, if $\sup_{l}|\mu_{\tau_{k_l}}'|(t) < +\infty$ for some subsequence $(\tau_{k_l})_{l\in\xN}, \ \tau_{k_l}\downarrow 0,$ and $t>0$, then $(L_F(u_{\tau_{k_l}}(t)))_{l\in\xN}$ is bounded in $\rm{W}^{1,1}(\Omega)$. In this case, by Rellich-Kondrachov Theorem, there exists a subsequence of $(L_F(u_{\tau_{k_l}}(t)))_{l\in\xN}$ which converges strongly in $\xLone(\Omega)$. If $\Big(L_F(u_{\tau_{k_{l_j}}}(t))\Big)_{j\in\xN}$ converges to some $\mathfrak{L}\in\xLone(\Omega)$, then it will, in turn, contain a subsequence which converges to $\mathfrak{L}$ pointwise $\mathscr{L}^d$-a.e., and by \eqref{eq: inverse LF}, the corresponding subsequence of $(u_{\tau_{k_{l_j}}}(t))_{j\in\xN}$ will converge to some $\mathfrak{u}: \Omega\to [0, +\infty)$ pointwise $\mathscr{L}^d$-a.e. and $\mathfrak{L} = L_F(\mathfrak{u})$. Using Egorov Theorem and the facts that $u_{\tau_k}(t)\stackrel{\xLone}{\rightharpoonup} u(t)$ and $\mathfrak{u}\in\xLone(\Omega)$ by Fatou's lemma, we obtain $\mathfrak{u}=u(t)$ and thus $\mathfrak{L}=L_F(u(t))$. This shows that the whole sequence $(L_F(u_{\tau_{k_l}}(t)))_{l\in\xN}$ converges to $L_F(u(t))$ strongly in $\xLone(\Omega)$ whenever $\sup_l |\mu_{\tau_{k_l}}'|(t) < +\infty$. Furthermore, in this case, the corresponding sequence $(\nabla L_F(u_{\tau_{k_l}}(t)))_{l\in\xN}$ of weak gradients is bounded in $\xLone(\Omega; \xR^d)$ and equiintegrable because $(u_{\tau_{k_l}}(t))_{l\in\xN}$ is equiintegrable and 
\begin{eqnarray*}
\int_E{\Big|\nabla L_F(u_{\tau_{k_l}}(t,x))\Big|\xdif x} \ \leq \ \Big(\sup_{x\in\Omega}|\nabla V(x)| \Big(\int_\Omega{u_{\tau_{k_l}}(t,x)\xdif x}\Big)^{1/2} + \frac{1}{\sqrt{\Lambda}}|\mu_{\tau_{k_l}}'|(t)\Big)\Big(\int_E{u_{\tau_{k_l}}(t,x)\xdif x}\Big)^{1/2}
\end{eqnarray*}
for every Borel set $E\subset \Omega$, by the preceding estimations of $\int_\Omega{|w_{\tau_k}^n(x)|^2u_{\tau_k}^n(x)\xdif x}, \ \nabla L_F(u_{\tau_k}^n) = w_{\tau_k}^n u_{\tau_k}^n$ and H\"older's inequality.  
 %we have 
%\begin{eqnarray*}
%&&\Big|\int_\Omega{-L_F(u(t,x))\xtr \xDif v(x) \xdif x}\Big| \ = \ \lim_{l\to\infty} \Big|\int_\Omega{-L_F(u_{\tau_{k_l}}(t,x))\xtr \xDif v(x) \xdif x}\Big| \\
%&&\leq \Big(\sup_{x\in\Omega}|\nabla V(x)| \Big(\int_\Omega{u(t,x)\xdif x}\Big)^{1/2} + \frac{1}{\sqrt{\Lambda}}\sup_{l\in\xN} |\mu_{\tau_{k_l}}'|(t)\Big)\Big(\int_\Omega{|v(x)|^2u(t,x)\xdif x}\Big)^{1/2} \\
%&&\leq  \Big(\sup_{x\in\Omega}|\nabla V(x)| \Big(\int_\Omega{u(t,x)\xdif x}\Big)^{1/2} + \frac{1}{\sqrt{\Lambda}} \sup_{l\in\xN} |\mu_{\tau_{k_l}}'|(t)\Big)\Big(\int_\Omega{u(t,x)\xdif x}\Big)^{1/2} \ \sup_{x\in\Omega}|v(x)| 
%\end{eqnarray*}   
%for every continuously differentiable function $v: \Omega\to\xR^d$ with compact support in $\Omega$, by the preceding estimations of $\Big|\int_\Omega{-L_F(u_{\tau_k}^n)\xtr \xDif v(x)\xdif x}\Big|$ and a corresponding passage to the limit $\tau_{k_l}\downarrow 0$; we may repeat the arguments showing that $L_F(u_{\tau_k}^n)\in\rm{W}^{1,1}(\Omega)$ in order to obtain
Dunford-Pettis Theorem and the above considerations show that 
$L_F(u(t))\in\rm{W}^{1,1}(\Omega)$ for every $t>0$ for which there exists a bounded subsequence of $(|\mu_{\tau_k}'|(t))_{k\in\xN}$ and $(\nabla L_F(u_{\tau_{k_l}}(t)))_{l\in\xN}$ converges to $\nabla L_F(u(t))$ weakly in $\xLone(\Omega; \xR^d)$ whenever $\sup_l |\mu_{\tau_{k_l}}'|(t) < +\infty$. By Fatou's lemma, this is true for a.e. $t>0$ because \eqref{eq: standard tools 1}, \eqref{eq: standard tools 2} yield 
\begin{equation}\label{eq: discrete metric derivative}
\Delta_T \ := \ \sup_{k\in\xN} \int_0^T{|\mu_{\tau_k}'|(t)^2\xdif t} \ < \ +\infty \quad\quad \text{ for all } T>0. 
\end{equation}    
Now, let $t>0$ and $\tau_{k_l}\downarrow 0$ such that $\sup_{l} |\mu_{\tau_{k_l}}'|(t) < +\infty$. Using Dunford-Pettis-Theorem, the above estimation of $\int_E{|\hat{L}_F(u_{\tau_k}^n(x))|\xdif x}$ and the equiintegrability of $(u_{\tau_{k_l}}(t))_{l\in\xN}$, we see that there exists a subsequence of $(\hat{L}_F(u_{\tau_{k_l}}(t)))_{l\in\xN}$ which converges weakly in $\xLone(\Omega)$. The preceding considerations show that every subsequence of $(u_{\tau_{k_l}})_{l\in\xN}$ contains a subsubsequence which converges to $u(t)$ pointwise $\mathscr{L}^d$-a.e.. We may use Egorov Theorem and the continuity of $\hat{L}_F$ in order to conclude that the whole sequence $(\hat{L}_F(u_{\tau_{k_l}}(t)))_{l\in\xN}$ converges to $\hat{L}_F(u(t))$ weakly in $\xLone(\Omega)$.
It is apparent from \eqref{eq: discrete metric derivative} and the preceding convergence results and estimations that $L_F(u)\in\xLtwo_{\mathrm{loc}}([0, +\infty); \rm{W}^{1,1}(\Omega))$ and $\hat{L}_F(u)\in\xLtwo_{\mathrm{loc}}([0, +\infty); \xLone(\Omega))$.  

All in all, we obtain 
\begin{eqnarray*}
&& \liminf_{k\to\infty} \Big[\mathfrak{I}_{\tau_k, \psi}(t) + \epsilon |\mu_{\tau_k}'|(t)^2 \Big] \\ && \geq \ \int_\Omega{\Big[\Lambda \Big\langle\nabla L_F(u(t,x)) + u(t,x)\nabla V(x),\nabla_x\psi(t,x)\Big\rangle + \Sigma\Big(\hat{L}_F(u(t,x)) + V(x)u(t,x)\Big) \psi(t,x)\Big]\xdif x}
\end{eqnarray*}
for every $\psi\in\xCtwo_{\mathrm{c}}(\xR\times\Omega), \ \epsilon > 0$ and almost every $t>0$, where $\mathfrak{I}_{\tau_k, \psi}(t)$ is defined as above, i.e. as
\begin{equation*}
\int_\Omega{\Big[-\Lambda L_F(u_{\tau_k}(t,x))\Delta_x\psi(t,x) + \Sigma\hat{L}_F(u_{\tau_k}(t,x))\psi(t,x) + \Big(\Lambda\langle \nabla V(x), \nabla_x\psi(t,x)\rangle + \Sigma V(x) \psi(t,x)\Big)u_{\tau_k}(t,x)\Big]\xdif x}. 
\end{equation*}   
Note that, by the above estimations of $\int_\Omega{|\hat{L}_F(u_{\tau_k}^n(x))|\xdif x}$ and $\Big|\int_\Omega{-L_F(u_{\tau_k}^n(x))\xtr \xDif v(x)\xdif x}\Big|$ (here we set $v(x)$ equal to $\nabla_x\psi(t,x)$) and by Cauchy's inequality with $\epsilon > 0$, we have 
\begin{eqnarray*}
&& \mathfrak{I}_{\tau_k, \psi}(t) + \epsilon |\mu_{\tau_k}'|(t)^2 \\
&& \geq -\Lambda\Big(C_V C_{T_\psi}^{1/2} + \frac{1}{\sqrt{\Lambda}}|\mu_{\tau_k}'|(t)\Big)C_{T_\psi}^{1/2}C_\psi - \Sigma C_\psi\Big(C_V C_{T_\psi} + \frac{1}{\sqrt{\Sigma}}|\mu_{\tau_k}'|(t)C_{T_\psi}^{1/2}\Big) - (\Lambda + \Sigma)C_V C_\psi C_{T_\psi} + \epsilon |\mu_{\tau_k}'|(t)^2 \\
&& \geq -2(\Lambda + \Sigma) C_V C_{T_\psi} C_\psi - \sqrt{\Lambda}|\mu_{\tau_k}'|(t)C_{T_\psi}^{1/2}C_\psi - \sqrt{\Sigma}|\mu_{\tau_k}|'(t)C_{T_\psi}^{1/2}C_\psi +\epsilon|\mu_{\tau_k}'|(t)^2 \\
&& \geq -(\Lambda + \Sigma)\Big(2C_V C_{T_\psi} C_\psi + \frac{1}{\epsilon} C_{T_\psi} C_\psi^2\Big) + \frac{\epsilon}{2}|\mu_{\tau_k}'|(t)^2 
\end{eqnarray*}
with $C_V:= \sup_{x\in\Omega}(|V(x)|+|\nabla V(x)|), \ C_\psi:=\sup_{(t,x)\in\xR\times\Omega}(|\psi(t,x)|+|\nabla_x \psi(t,x)|), \ T_\psi > 0$ such that $\psi(t, \cdot)\equiv 0$ for all $t\geq T_\psi$, and $C_{T_\psi}:=\sup\Big\{\int_\Omega{u_{\tau_k}(t,x)\xdif x }: \ k\in\xN, \ t\in (0, T_\psi)\Big\}$ (which is finite by \eqref{eq: standard tools 1}), so that the limit inferior of $\Big(\mathfrak{I}_{\tau_k, \psi}(t) + \epsilon |\mu_{\tau_k}'|(t)^2 \Big)_{k\in\xN}$ is indeed either $+\infty$ (if there is no bounded subsequence of $(|\mu_{\tau_k}'|(t))_{k\in\xN}$) or the limit of some subsequence for which the corresponding subsequence of $(|\mu_{\tau_k}'|(t))_{k\in\xN}$ is bounded. We may apply Fatou's lemma in order to obtain
\begin{eqnarray*}
\liminf_{k\to\infty} \int_0^\infty{\mathfrak{I}_{\tau_k, \psi}(t)\xdif t}  +  \epsilon\Delta_{T_\psi}  \geq  \liminf_{k\to\infty} \int_0^{T_\psi}{\Big[\mathfrak{I}_{\tau_k, \psi}(t) + \epsilon |\mu_{\tau_k}'|(t)^2\Big]\xdif t}  \geq  \int_0^{T_\psi}{\liminf_{k\to\infty}\Big[\mathfrak{I}_{\tau_k, \psi}(t) + \epsilon |\mu_{\tau_k}'|(t)^2\Big]\xdif t} \\
\geq \int_0^{T_\psi}{\int_\Omega{\Big[\Lambda \Big\langle\nabla L_F(u(t,x)) + u(t,x)\nabla V(x),\nabla_x\psi(t,x)\Big\rangle + \Sigma\Big(\hat{L}_F(u(t,x)) + V(x)u(t,x)\Big) \psi(t,x)\Big]\xdif x}\xdif t}  =  \mathcal{I}_{F,V,\psi, u}
\end{eqnarray*} 
with $\mathcal{I}_{F,V,\psi, u}$ defined as in \eqref{eq: weak form 2} and $\Delta_{T_\psi}$ according to \eqref{eq: discrete metric derivative}. This shows that
\begin{equation*}
\lim_{k\to\infty} \int_0^\infty{\mathfrak{I}_{\tau_k, \psi}(t)\xdif t} \ = \ \mathcal{I}_{F,V,\psi, u} \quad \text{ for every } \psi\in \xCtwo_{\mathrm{c}}(\xR\times\Omega)
\end{equation*}
because we may let $\epsilon\downarrow 0$ and switch between $\psi$ and $-\psi$. The proof of \eqref{eq: weak form 1} for $\psi\in\xCtwo_{\mathrm{c}}(\xR\times \Omega)$ is complete.   \\

Now, we include the no-flux boundary condition \eqref{eq: our bc} in a weak form and prove \eqref{eq: weak form 1} for all $\psi\in\xCtwo_{\mathrm{c}}(\xR\times\xR^d)$. The reason why we cannot just repeat the previous proof for $\psi\in\xCtwo_{\mathrm{c}}(\xR\times\Omega)$ is that the derivation of \eqref{eq: cond of first order} from Prop. \ref{prop: subdif of HK}, Ass. \ref{ass: F} and \eqref{eq: diff cond V} will fail for general continuously differentiable functions $v$ with compact support in $\xR^d$ (but not in $\Omega$), cf. Ex. \ref{ex: bc} below. However, using the preceding results, a necessary condition of first order can still be obtained in this case. For $k\in\xN, \ n\in\xN,$ let $\gamma_{\tau_k}^n\in\mathcal{M}(\bar{\Omega}\times \bar{\Omega})$ be optimal in the definition of $\HK_{\Lambda, \Sigma}(\mu_{\tau_k}^n, \mu_{\tau_k}^{n-1})^2$ according to \eqref{eq: LET}, with first marginal $\gamma_{\tau_k, 1}^n\ll\mu_{\tau_k}^n$ and second marginal $\gamma_{\tau_k, 2}^n\ll\mu_{\tau_k}^{n-1}$, and Lebesgue decompositions 
\begin{equation*} 
\mu_{\tau_k}^n = \rho_{\tau_k, 1}^n\gamma_{\tau_k, 1}^n + (\mu_{\tau_k}^n)^\bot \quad\text{ and } \quad \mu_{\tau_k}^{n-1}=\rho_{\tau_k, 2}^{n}\gamma_{\tau_k, 2}^{n}+(\mu_{\tau_k}^{n-1})^\bot,
\end{equation*}
cf. \eqref{eq: Lebesgue decomposition}.  
By Cor. \ref{cor: subdif of HK}, Ass. \ref{ass: F}, \eqref{eq: diff cond V}, the same arguments showing \eqref{eq: cond of first order} and as $L_F(u_{\tau_k}^n)\in\rm{W}^{1,1}(\Omega)$, we obtain the following necessary condition of first order 
 \begin{equation}\label{eq: cond of first order bc}
\int_{\Omega}{\Big[\langle \nabla L_F(u_{\tau_k}^n), v \rangle + 2\hat{L}_F(u_{\tau_k}^n)R + \Big(\langle\nabla V, v\rangle + 2VR\Big)u_{\tau_k}^n\Big]\xdif x} \ = \ \frac{1}{\tau_k} \Big( \mathfrak{F}_{\tau_k, n, v, R} \ - \ \frac{4}{\Sigma}\int_{\bar{\Omega}}{R(x)\xdif(\mu_{\tau_k}^n)^\bot}\Big)
\end{equation}  
for all $\xCinfty$-functions $v: \Omega\to\xR^d$ with compact support in $\Omega$ and all bounded Borel measurable functions $R: \Omega\to\xR$, where $\mathfrak{F}_{\tau_k, n, v, R}$ is defined as
\begin{eqnarray*} && \frac{4}{\Sigma}\int_{\bar{\Omega}\times \bar{\Omega}}{\Big[-\rho_{\tau_k,1}^n(x_1)R(x_1)+\sqrt{\rho_{\tau_k, 1}^n(x_1)\rho_{\tau_k, 2}^n(x_2)}R(x_1)\cos(\sqrt{\Sigma/4\Lambda}||x_1-x_2||)\Big]\xdif \gamma_{\tau_k}^n} \\  &+& \frac{4}{\Sigma}\int_{\bar{\Omega}\times\bar{\Omega}}{\sqrt{\Sigma/4\Lambda \ \rho_{\tau_k, 1}^n(x_1)\rho_{\tau_k, 2}^n(x_2)}\ \langle S_{\Lambda,\Sigma}(x_1,x_2), v(x_1)\rangle\xdif\gamma_{\tau_k}^n}.   
\end{eqnarray*} 
According to \eqref{eq: LET},  Thm. 4.5 in \cite{Gangbo-McCann96} and Thm. 6.6 in \cite{liero2018optimal}, there exist a Borel function $\sigma_{\tau_k, 1}^n: \bar{\Omega}\to[0, +\infty)$ and a Borel optimal transport mapping $t_{\tau_k}^n: \bar{\Omega}\to\bar{\Omega}$ such that
\begin{equation}\label{eq: optimal transport map}
\gamma_{\tau_k, 1}^n = \sigma_{\tau_k, 1}^n\mu_{\tau_k}^n, \quad\quad\quad \gamma_{\tau_k}^n \ = \ (I\times t_{\tau_k}^n)_{\#} \gamma_{\tau_k,1}^{n} \ = \ (I\times t_{\tau_k, 1}^n)_{\#}(\sigma_{\tau_k, 1}^n u_{\tau_k}^n \mathscr{L}^d). 
\end{equation}
Setting $R\equiv 0$ in \eqref{eq: cond of first order bc} and applying \eqref{eq: optimal transport map}, we obtain
\begin{eqnarray*}
&& \int_\Omega{\langle \nabla L_F(u_{\tau_k}^n(x)) + u_{\tau_k}^n(x)\nabla V(x), v(x)\rangle \xdif x} \\ && = \ \frac{1}{\tau_k} \frac{4}{\Sigma}\int_{\Omega}{\sqrt{\Sigma/4\Lambda \ \rho_{\tau_k, 1}^n(x)\rho_{\tau_k, 2}^n(t_{\tau_k}^n(x))}\ \langle S_{\Lambda,\Sigma}(x,t_{\tau_k}^n(x)), v(x)\rangle \sigma_{\tau_k,1}^n(x) u_{\tau_k}^n(x)\xdif x} 
\end{eqnarray*} 
for all $\xCinfty$-functions $v: \Omega\to\xR^d$ with compact support in $\Omega$. This shows that
\begin{equation*}
\nabla L_F(u_{\tau_k}^n(x)) + u_{\tau_k}^n(x)\nabla V(x) \ = \ \frac{2}{\tau_k\sqrt{\Sigma\Lambda}} \sqrt{\rho_{\tau_k, 1}^n(x)\rho_{\tau_k, 2}^n(t_{\tau_k}^n(x))} \ \sigma_{\tau_k,1}^n(x) u_{\tau_k}^n(x) S_{\Lambda,\Sigma}(x,t_{\tau_k}^n(x)) \quad \mathscr{L}^d\text{-a.e. in } \Omega.
\end{equation*}
Consequently, \eqref{eq: cond of first order bc} holds good for all bounded Borel measurable functions $v: \Omega\to\xR^d$ and $R:\Omega\to\xR$. Applying the proof of Cor. \ref{cor: subdif of HK}, Prop. \ref{prop: connection of subdif with MM} and the preceding arguments for $\psi\in\xCtwo_{\mathrm{c}}(\xR\times\Omega)$, we obtain \eqref{eq: weak form 1} for all $\psi\in\xCtwo_{\mathrm{c}}(\xR\times\xR^d)$. The proof of Thm. \ref{thm: main thm} is complete. 
\end{proof}

The following example shows that Ass. \ref{ass: F} does \textit{not} imply such formula \eqref{eq: diff cond} for every continuously differentiable function $v: \xR^d\to\xR^d$ with compact support in $\xR^d$. 
\begin{xmpl}\label{ex: bc}
We can identify $\nu_0 = u_0\mathscr{L}^d\in\{\mathcal{F} < +\infty\}$ with a nonnegative finite Radon measure $\nu_0\in\mathcal{M}(\xR^d)$ by setting $u_0\equiv 0$ outside $\Omega$ and define the curve $h\mapsto \nu_h\in\mathcal{M}(\xR^d)$ as in Ass. \ref{ass: F}, for a continuously differentiable function $v: \xR^d\to\xR^d$ with compact support in $\xR^d$ and a bounded Borel measurable function $R: \xR^d\to\xR^d$. Again, we obtain \eqref{eq: 78} which holds good on $\xR^d$ for now. The measures $\nu_h$ can be restricted to measures $\nu_h\in\mathcal{M}(\bar{\Omega})$ and $\mathcal{F}(\nu_h) := \int_\Omega{F(u_h(x))\xdif x}$. 

Now, let $\Omega = (0,1), \ F(s):=s^2$ and $v: \xR\to\xR$ be a continuously differentiable function satisfying
\begin{equation*}
v(x) \begin{cases}
= 0 &\text{ if } |x|\geq 2, \\
\in [-1,0) &\text{ if } -2 < x < 1, \\
= -1 &\text{ if } 0\leq x < \frac{1}{2}, \\
\geq 0 &\text{ if } x\geq 1.  
\end{cases}
\end{equation*}
Then, for $h>0$ small enough, we have
\begin{eqnarray*}
\frac{\mathcal{F}(\nu_h) - \mathcal{F}(\nu_0)}{h} \ = \ -\frac{1}{h}\int_0^h{F(u_0(x))\xdif x} + \underbrace{\frac{1}{h}\int_h^1{\Big(F\Big(\frac{u_0(x)(1+hR(x))^2}{1+hv'(x)}\Big)(1+hv'(x))-F(u_0(x))\Big) \xdif x}}_{\to \ \int_\Omega{u_0(x)^2(4R(x)-v'(x))\xdif x} \text{ as } h\downarrow 0}
\end{eqnarray*}
and we note that the first term on the right-hand side cannot be controlled; if we take $u_0(x):= x^{-1/4}, \ x\in(0,1)$, then $\nu_0=u_0\mathscr{L}^1\in\{\mathcal{F} <+\infty\}$ and 
\begin{equation*}
-\frac{1}{h}\int_0^h{F(u_0(x))\xdif x} \ = \ -\frac{2\sqrt{h}}{h} \ \to \ -\infty \quad \text{ as } h\downarrow 0. 
\end{equation*}
\end{xmpl}

\subsection{Comments and outlook}\label{subsec: 3.3}
The discussion in Sect. \ref{subsec: 3.1} and Prop. \ref{prop: typical ass} show that, under the assumptions of Thm. \ref{thm: main thm}, every Generalized Minimizing Movement $\mu\in\mathrm{GMM}(\Phi; \mu_0), \ \mu_0\in\{\mathcal{E} < +\infty\}$, associated with $\Phi$ and $\mathcal{E}$ as in \eqref{eq: our MM} and \eqref{eq: our E}, is locally absolutely continuous and satisfies the energy dissipation inequality \eqref{eq: energy inequality}, i.e.
\begin{equation*}
\mathcal{E}(\mu_0) - \mathcal{E}(\mu(t)) \ \geq \ \frac{1}{2} \int_0^t{|\partial^-\mathcal{E}|(\mu(r))^2 \mathrm d r} + \frac{1}{2}\int_0^t{|\mu'|(r)^2 \mathrm d r}
\end{equation*}
for all $t > 0$. In this paper, we have restricted ourselves to proving that our Minimizing Movement scheme yields weak solutions to a class of reaction-diffusion equations. It will be worth studying the corresponding energy dissipation inequalities in a subsequent paper as they will provide additional information (cf. Ex. \ref{ex: 3.10} below).

\begin{dfntn}[Absolutely continuous curves, relaxed slope]
Let $(\mathscr{S}, d)$ be a complete metric space. We say that a curve $u: [0, +\infty) \to \mathscr{S}$ is locally absolutely continuous if there exists $m\in \xLone_{\mathrm{loc}}(0, +\infty)$ such that 
\begin{equation*}
d(u(s),u(t)) \leq \int^{t}_{s}{m(r) \xdif r} \quad \quad \text{ for all } 0 \leq s\leq t < +\infty. 
\end{equation*}\\
In this case, the limit
\begin{equation*}
|u'|(t) := \mathop{\lim}_{s\to t} \frac{d(u(s),u(t))}{|s-t|}
\end{equation*}
exists for $\mathscr{L}^1$-a.e. $t$, the function $t \mapsto |u'|(t)$ belongs to $\xLone_{\mathrm{loc}}(0, +\infty)$ and is called the metric derivative of $u$. The metric derivative is $\mathscr{L}^1$-a.e. the smallest admissible function $m$ in the definition above.

Let $\mathcal{E}: \mathscr{S} \to (-\infty, +\infty]$ be given. 
We define the local slope at $x\in \{\mathcal{E} < +\infty\}$ as
\begin{equation*}
|\partial\mathcal{E}|(x) := \mathop{\limsup}_{d(x,y)\to 0} \frac{(\mathcal{E}(x)-\mathcal{E}(y))^+}{d(x,y)} 
\end{equation*}
and the relaxed slope $|\partial^- \mathcal{E}|: \mathscr{S}\to [0, +\infty]$ of $\mathcal{E}$ as 
\begin{equation*}
|\partial^- \mathcal{E}|(x) := \inf  \left\{\mathop{\liminf}_{n\to \infty} |\partial\mathcal{E}|(x_n): \ d(x_n, x)\to 0, \ \sup_{n}\mathcal{E}(x_n) < +\infty \right\}.
\end{equation*}
\end{dfntn}

We refer to (\cite{AmbrosioGigliSavare05}, Chaps. 1 and 2) for a detailed account of these and further definitions which are important in connection with the characterization of gradient flows in metric spaces by such energy dissipation (in)equality (cf. introductory part). 

\begin{xmpl}\label{ex: 3.10}
Let $d\in\{1, 2\}$. For $F: [0, +\infty)\to\xR, \ F(s):= -\sqrt{s}+s^p \ (p > 1)$ and $V\equiv 0$, we define $\mathcal{E}$ and $\Phi$ as in \eqref{eq: our E} and \eqref{eq: our MM}. In this case, the functional $\mathcal{E}$ is geodesically convex on $(\mathcal{M}(\bar{\Omega}), \HK_{\Lambda, \Sigma})$ according to \cite{mlsHKconvex} and 
\begin{eqnarray*}
|\partial^-\mathcal{E}|(\mu) = \begin{cases} |\partial\mathcal{E}|(\mu) &\text{ if } \mu\in\{\mathcal{E} < +\infty\}, \\
+\infty &\text{ else }
\end{cases}
\end{eqnarray*}
by Cor. 2.4.10 in \cite{AmbrosioGigliSavare05}. 
Obviously, $u\equiv 0$ is a solution to \eqref{eq: weak form 1}, \eqref{eq: weak form 2}. However, Thm. \ref{thm: main thm} (which is applicable in this example, cf. Ex. \ref{ex: 3.8}) will not yield this trivial solution, i.e. 
\begin{equation*}
\eta_0\notin\mathrm{GMM}(\Phi; \eta_0) 
\end{equation*}
(where $\eta_0$ denotes the null measure). Indeed, setting $\eta_N:=u_N\mathscr{L}^d, \ u_N\equiv \frac{1}{N}$ we can easily compute that
\begin{equation*}
|\partial^-\mathcal{E}|(\eta_0) = |\partial\mathcal{E}|(\eta_0) \geq \lim_{N\to\infty}\frac{-\mathcal{E}(\eta_N)}{\HK_{\Lambda, \Sigma}(\eta_N, \eta_0)} = \lim_{N\to\infty}\frac{(1/\sqrt{N} - 1/N^p)\mathscr{L}^d(\Omega)}{2/\sqrt{\Sigma} \sqrt{1/N \ \mathscr{L}^d(\Omega)}} = \frac{\sqrt{\Sigma \mathscr{L}^d(\Omega)}}{2} \ > \ 0,
\end{equation*} 
which shows that the constant curve $\eta(t)\equiv \eta_0$ does not satisfy the corresponding energy dissipation inequality.     
\end{xmpl}  

Our next comment concerns the initial data. We may replace $\mu_0=u_0\mathscr{L}^d\in\{\mathcal{E} < +\infty\}$ in the Minimizing Movement scheme \eqref{eq: MM scheme} associated with \eqref{eq: our MM}, \eqref{eq: our E} by a sequence $(\mu_\tau^0)_\tau$ of measures $\mu_\tau^0=u_\tau^0\mathscr{L}^d$ satisfying $u_\tau^0\stackrel{\xLone}{\rightharpoonup} u_0, \ \sup_\tau \mathcal{F}(\mu_\tau^0) < +\infty$, and still obtain the same results as in Thm. \ref{thm: main thm}. 

We remark that we have left aside the possibility of adding an interaction energy functional to $\mathcal{E}$ for the sake of clear presentation. 

We expect that our arguments will form the basis for a Minimizing Movement approach to scalar reaction-diffusion equations in other settings, too, e.g. if $X=\xR^d$ or $X$ is a subset of a general separable Hilbert space or the energy functional is modified. We do not want to expound on how to adapt our assumptions and our proof for such cases, just give an example of suitable assumptions if $X=\xR^d$. 
\begin{xmpl}
We suppose that $F: [0, +\infty) \to \xR$ is continuous, strictly convex, differentiable in $(0, +\infty)$, has superlinear growth \eqref{eq: superlinear growth} and satisfies $F(0)=0, \ F(s)\geq -C_F s$ (for some $C_F>0$). Let $V: \xR^d\to\xR$ be locally Lipschitz continuous and let us suppose that $V\geq 0$ and $V(x)\to +\infty$ if $|x|\to+\infty$. We define $\mathcal{E}:=\mathcal{F}+\mathcal{V}: \mathcal{M}(\xR^d)\to (-\infty, +\infty]$ and $\Phi$ as in \eqref{eq: our E} and \eqref{eq: our MM} with $\Omega$ replaced by $\xR^d$, and we suppose that $\mathcal{F}$ satisfies a differentiability assumption which is like Ass. \ref{ass: F} (but with $\Omega$ and $\xLone(\Omega)$ replaced by $\xR^d$ and $\xLone_{\mathrm{loc}}(\xR^d)$ respectively, and for continuously differentiable functions $v: \xR^d\to\xR^d$ and bounded Borel functions $R:\xR^d\to\xR$, $v$ and $R$ both with compact support in $\xR^d$). 
Then similar arguments as in Sect. \ref{subsec: 3.2} will show that the associated Minimizing Movement approach yields weak solutions to the corresponding scalar reaction-diffusion equation on $\xR^d$; the results are similar to those of Thm. \ref{thm: main thm}. %(with $\Omega$ replaced by $\xR^d$ and $L_F(u(t))\in\rm{W}_{\mathrm{loc}}^{1,1}(\xR^d)$ for a.e. $t>0$). 
Note that the growth condition on $V$ makes an application of Dunford-Pettis-Theorem on $\xR^d$ possible.  
\end{xmpl}      
%Also, the energy dissipation inequality can still be obtained in cases in which our proof of the reaction-diffusion equation (in the weak form \eqref{eq: weak form 1}, \eqref{eq: weak form 2}) fails, e.g. if condition \eqref{eq: technical ass} does not hold. 
\begin{acknowledgement}
\paragraph{\textit{Acknowledgements}}
Giuseppe Savar\'e suggested this topic to me and I had the opportunity to discuss parts of it with him, Alexander Mielke and Martin Brokate; I would like to express my thanks to them.  
I gratefully acknowledge support from the Erwin Schr\"odinger International Institute for Mathematics and Physics (Vienna) during my participation in the programme ``Optimal Transport''.
\end{acknowledgement}
\bibliographystyle{siam}
\bibliography{bibliografia}
\end{document}